\documentclass[10pt]{amsart}

\usepackage{amssymb}
\usepackage{amscd}
\usepackage{upref}
\usepackage{tikz}
\usepackage[all]{xy}

\theoremstyle{plain}
\newtheorem{theorem}{Theorem}[section]
\newtheorem{corollary}[theorem]{Corollary}
\newtheorem{lemma}[theorem]{Lemma}
\newtheorem{proposition}[theorem]{Proposition}

\newtheorem*{Theorem A}{Theorem A}
\newtheorem*{Theorem B}{Theorem B}
\newtheorem*{Theorem C}{Theorem C}
\newtheorem*{Vanishing Theorem}{Vanishing Theorem}
\newtheorem*{Harer's Stability Theorem}{Harer's Stability Theorem}

\theoremstyle{definition}

\theoremstyle{remark}

\newtheorem{remark}[theorem]{Remark}

\numberwithin{equation}{section}

\newcommand{\HC}{\mathcal H_*\mathcal C}
\newcommand{\Hc}{\mathcal H_0\mathcal C}

\begin{document}

\title[Stable String Operations Are Trivial]{Stable  String Operations
are Trivial}

\author{Hirotaka Tamanoi }

\address{ Department of Mathematics,
 University of California Santa Cruz  \newline
\indent Santa Cruz, CA 95064}
\email[]{tamanoi@math.ucsc.edu}
\date{}
\begin{abstract}
We show that in closed string topology and in open-closed string topology with one $D$-brane, higher genus stable string operations are trivial. This is a consequence of Harer's stability theorem and related stability results on the homology of mapping class groups of surfaces with boundaries. In fact, this vanishing result is a special case of a general result which applies to all homological conformal field theories with a property that in the associated topological quantum field theories, the string operations associated to genus one cobordisms with one or two boundaries vanish. In closed string topology, the base manifold can be either finite dimensional, or infinite dimensional with finite dimensional cohomology for its based loop space. The above vanishing result is based on the triviality of string operations associated to homology classes of mapping class groups which are in the image of stabilizing maps. 
\end{abstract}

\maketitle

\tableofcontents

\section{Introduction}

Let $M$ be a closed oriented smooth $d$ dimensional manifold  and let $LM$ be the loop space consisting of continuous maps from $S^1$ into $M$. In this paper, we use (co)homology with integral coefficients unless otherwise stated, except in section 4. Chas and Sullivan \cite{CS} showed that the homology of the loop space with a degree shift $\mathbb H_*(LM)=H_{*+d}(LM)$ has the structure of a Batalin-Vilkovisky algebra. Cohen and Godin \cite{CG}, Godin \cite{Go2}, and Cohen and Schwarz \cite{CoSch} showed that $\mathbb H_*(LM)$ even carries the structure of a homological conformal field theory (HCFT). 

Namely, let $F_{g,p+q}$ be a connected smooth oriented genus $g$ surface with $p+q$ parametrized boundaries out of which $p$ of them are designated as incoming and the other $q$ as outgoing. The mapping class group $\Gamma_{g,p+q}$ is the group of isotopy classes of orientation preserving diffeomorphisms of $F_{g,p+q}$ fixing boundaries pointwise. We assume $p\ge0, q\ge1$. This condition is called the positive boundary condition. Then to each homology class of the mapping class group, HCFT structure assigns the following string operation acting on $H_*(LM)$ (modulo K\"unneth theorem): 
\begin{equation}\label{string operation}
\mu_{g,p+q}: H_*(B\Gamma_{g,p+q})\otimes H_*(LM)^{\otimes p} \longrightarrow  H_*(LM)^{\otimes q}, 
\end{equation} 
lowering degree by $-d\cdot \chi(F_{g,p+q})=d(2g+p+q-2)$, where $B\Gamma_{g,p+q}$ is the classifying space of the mapping class group $\Gamma_{g,p+q}$. Since we assume $q\ge1$, $B\Gamma_{g,p+q}$ is homotopy equivalent to the moduli space $\mathfrak M_{g,p+q}$ of connected Riemann surfaces of genus $g$ with $p+q$ disjoint holomorphically embedded discs. A representation theory of these moduli spaces $\mathfrak M_{g,p+q}$ is a conformal field theory. Hence the name homological conformal field theory is a suitable one in our framework. 

When the manifold $M$ is infinite dimensional, Cohen-Godin's construction of string operations does not work. However, if $M$ is simply connected and the cohomology of its based loop space $\Omega M$ is finite dimensional over a field $k$, then Chataur and Menichi \cite{CM} constructed a homological conformal field theory structure on the cohomology $H^*(LM;k)$ under a different condition $p,q\ge1$ (noncompact HCFT). String operations in this case, are of the following form:
\begin{equation}\label{cohomology string operation}
\mu_{g,p+q}: H_*(B\Gamma_{g,p+q};k)\otimes H^*(LM;k)^{\otimes q} \longrightarrow  H^*(LM;k)^{\otimes p}.
\end{equation}
In section \ref{vanishing theorem (II)}, we show that the coproduct for $H^*(LM;k)$ is trivial (i) when $H^*(\Omega M;\mathbb Z)$ is torsion free and $k$ is any field, and (ii) when $H^*(\Omega M;\mathbb Z)$ is $p$-torsion free where $p$ is the characteristic of the field $k$. See Corollary \ref{trivial coproduct}. The homology $H_*(LM;k)$ has a dual HCFT structure with trivial product under the same condition on $\Omega M$. Their construction applies in particular to classifying spaces $BG$ of finite dimensional Lie groups $G$. 
 
The main result of this paper is that all ``stable" string operations vanish both in closed string topology and in open-closed string topology. To describe this stability condition, we recall Harer's stability theorem \cite{H}. Let $T$ be a torus with two boundaries. By sewing one boundary of $T$ to any boundary of $F_{g,p+q}$, we obtain a surface $F_{g+1,p+q}$ of genus $g+1$. By extending a self-diffeomorphism of $F_{g,p+q}$ to one of $F_{g+1,p+q}$ by letting the extension to be  identity on $T$, we get an inclusion of groups $\varphi:\Gamma_{g,p+q} \rightarrow  \Gamma_{g+1,p+q}$. The induced map in homology is an isomorphism in low degrees, increasing with genus. Thus, the homology of mapping class groups stabilizes with increasing genus. More precisely, 

\begin{Harer's Stability Theorem}[Harer \cite{H} and Ivanov \cite{I1}, \cite{I2}]  Assume $p+q\ge1$. The stabilizing homomorphism 
\begin{equation*}
\varphi_*: H_k(B\Gamma_{g,p+q}) \longrightarrow H_k(B\Gamma_{g+1,p+q})
\end{equation*}
is an isomorphism and independent of $p+q$ for $g\ge2k+1$. It is onto and independent of $p+q$ for $g\ge2k$. 
\end{Harer's Stability Theorem} 

Thus for sufficiently large $g$, the homology of mapping class groups is independent of genus and the number of boundaries $p+q\ge1$. We will show that there are actually $p+q$ possibly different choices of stabilizing maps related by $\Sigma_p\times\Sigma_q$ equivariance, corresponding to different choices of boundary circles of $F_{g,p+q}$ used for sewing with $T$. The above stability theorem is valid for each of these maps. In fact, it turns out that in stable range, all of these $p+q$ choices give the same stabilizing map. This is a consequence of Harer's stability theorem. See Remark \ref{trivial S_n action}. 

For the statement of Ivanov's reformulation of homology stability in \cite{I2}, see Theorem \ref{Ivanov stability}. This reformulation is more convenient for us when we discuss open-closed string topology operations. The above stability range was proved by Ivanov improving over the original Harer's stability range. Note that $\varphi_*$ is always isomorphism for $k=0$ because all of the spaces $B\Gamma_{g,p+q}$ are connected.

We say that the homology group $H_k(B\Gamma_{g,p+q})$ is in stable range if the stabilizing map $\varphi_*$ into $H_k(B\Gamma_{g,p+q})$ is surjective, and all the subsequent $\varphi_*$ are isomorphisms, as in the following sequence: 
\begin{equation*}
H_k(B\Gamma_{g-1,p+q}) \xrightarrow[\text{onto}]{\varphi_*}
H_k(B\Gamma_{g,p+q}) \xrightarrow[\cong]{\varphi_*} H_k(B\Gamma_{g+1,p+q}) \xrightarrow[\cong]{\varphi_*} \dots. 
\end{equation*}
By the Ivanov's result, $H_k(B\Gamma_{g,p+q})$ is in stable range when $g\ge 2k+1$ for all $k\ge0$. 

\begin{Vanishing Theorem}[\textbf{Closed String Topology Case}] 
\textup{(I)} Let $M$ be a finite dimensional closed smooth oriented manifold. Consider the closed string topology for $M$. 

\textup{(i)} String operations \eqref{string operation} on $H_*(LM)$ associated to elements in the image $\textup{Im}\,\varphi_*\subset H_k(B\Gamma_{g,p+q})$ of any stabilizing map $\varphi_*$ are trivial. 

\textup{(ii)} String operations associated to any elements in the homology  $H_k(B\Gamma_{g,p+q})$ in stable range are all trivial. 

\noindent\textup{(II)} Let $M$ be a simply connected infinite dimensional space such that $H^*(\Omega M;k)$ is finite dimensional for some field $k$. Then with respect to \textup{(}co\textup{)}homology with coefficients in $k$, the above statements \textup{(i)} and \textup{(ii)} for string operations associated to elements in $H_k(B\Gamma_{g,p+q})$ acting on $H^*(LM;k)$ are valid. 
\end{Vanishing Theorem}

The second parts (ii) describe the meaning of the title of this paper. Obviously parts (ii) are consequences of parts (i) in view of Harer's Stability Theorem. Parts (i) apply to all stable operations as well as the majority of higher genus unstable operations, and consequently,  most higher genus string operations vanish. In the last section of this paper, we compute some genus one unstable string operations and show them to be trivial. 

In \cite{Go2}, higher string operations are also constructed in open-closed string topology for a finite dimensional manifold $M$ in which the set of $D$-branes (a collection of submanifolds of $M$ in which open strings can end) consists of just $M$. To describe these operations, let $S$ be a connected open-closed cobordism of genus $g(S)$ with $p$ incoming closed strings, $q$ outgoing closed strings, $r$ incoming open strings, $s$ outgoing open strings, and $m$ completely free boundaries. Let $\Gamma(S)$ be the mapping class group of $S$, where  diffeomorphisms are allowed to permute completely free boundaries carrying the same label.  Let $\sigma_m\cong\mathbb Z$ be the sign representation of the symmetric group $\Sigma_m$ on $m$ letters. Since there is a canonical surjective map $\Gamma(S)\to \Sigma_m$, the module $\sigma_m$ is also a module over $\Gamma(S)$. See section \ref{open-closed string topology} for details. Then the open-closed string operation is of the following form (modulo K\"unneth theorem):
\begin{equation}\label{open-closed operation}
\mu: H_*(\Gamma(S);\sigma_m^d)\otimes H_*(LM)^{\otimes p}\otimes H_*(M)^{\otimes r} \longrightarrow 
H_*(LM)^{\otimes q}\otimes H_*(M)^{\otimes s},
\end{equation}
where $d=\dim M$ and $\sigma_m^d=(\sigma_m)^{\otimes d}$. We prove in Proposition \ref{chi_S} that $(\det\chi_S)^d$ appearing in \cite{Go2} is the same as $\sigma_m^d$ as $\Gamma(S)$-modules. We will formulate and prove a stability property of the homology of the mapping class group $H_*(\Gamma(S);\sigma_m^r)$ with coefficients in $\sigma_m^r$ for $r\ge0$. In particular, we show that the group $H_k(\Gamma(S);\sigma_m^d)$ is in stable range when $g(S)\ge2k+1$. See Remark \ref{related stability} for related works. 

\begin{Vanishing Theorem}[\textbf{Open-Closed String Topology Case}] Let $M$ be a $d$-dimensional closed oriented smooth manifold. Consider the open-closed string topology for $M$ with $D$-brane set consisting only of $M$.  

\textup{(i)} Open-closed string operations \eqref{open-closed operation} associated to elements in $H_k(\Gamma(S);\sigma_m^d)$ in the image of any stabilizing map $\varphi_*$ are trivial. 

\textup{(ii)} Open closed string operations associated to any elements in the homology group $H_k(\Gamma(S);\sigma_m^d)$ in stable range are trivial. 
\end{Vanishing Theorem} 
At the end of section \ref{open-closed string topology}, we will comment on the general open-closed string topology case with an arbitrary collection of $D$-brane submanifolds. 

Let $\Gamma_{\infty,r}=\lim_{g\to\infty}\Gamma_{g,r}$ be the stable mapping class group for $r\ge1$. In view of the Harer's stability theorem, the homology of $\Gamma_{\infty,r}$ is independent of $r\ge1$. The stable homology of the mapping class groups is given by the homology $H_*(B\Gamma_{\infty,r})$ of the stable mapping class group. The homotopy type of $B\Gamma_{\infty}=B\Gamma_{\infty,r}$ has been identified by Madsen and Weiss \cite{MW}: they showed that $H_*(\mathbb Z\times B\Gamma_{\infty};\mathbb Z)\cong H_*(\Omega^{\infty}\mathbb CP^{\infty}_{-1};\mathbb Z)$, and as a consequence, its rational cohomology is given by
\begin{equation*}
H^*(B\Gamma_{\infty};\mathbb Q)\cong\mathbb Q[\kappa_1,\kappa_2,\dots,\kappa_n,\dots],
\end{equation*}
solving Mumford conjecture, where $\kappa_n$'s of degree $2n$ classes are Miller-Morita-Mumford classes.  The mod $p$ homology of $B\Gamma_{\infty}$ was computed by Galatius \cite{Ga}. In contrast, homology of mapping class groups in unstable range has not been well understood, and only a few groups have been calculated. Although statements in part (ii) of the vanishing theorems are not the same as saying that string operations associated to $H_*(B\Gamma_{\infty})$ are trivial, this would be a concise statement. 

The organization of this paper is as follows. In section 2, we recall homological conformal field theory and discuss some of its properties. In section 3, we analyze the meaning of Harer's Stability Theorem in a general homological conformal field theoretic context. From this, the vanishing theorem of string operations for finite dimensional manifolds follows as a special case of a general fact. In section \ref{vanishing theorem (II)}, we prove the vanishing theorem for infinite dimensional manifolds $M$ with finite dimensional $H^*(\Omega M;k)$. This is done by showing that the genus $1$ TQFT operator is trivial (Theorem \ref{vanishing of genus one map}). We also show that the coproduct in the loop cohomology $H^*(LM;k)$ is trivial and the Serre spectral sequence for the fibration $p:LM \to M$ collapses when $H^*(\Omega M;k)$ is an exterior algebra (Theorem \ref{cohomology of LM}). In particular, the coproduct in $H^*(LM;k)$ is trivial if $H^*(\Omega M;\mathbb Z)$ is $p$ torsion free, where $p$ is the characteristic of the coefficient field $k$ (Corollary \ref{trivial coproduct}). In section \ref{open-closed string topology}, we prove the corresponding vanishing theorem in open-closed string topology with a single $D$-brane consisting of $M$ itself, which follows from a stability result of certain mapping class groups with nontrivial module coefficients obtained by a spectral sequence comparison argument. In the last section, we compute some genus one unstable string operations, and show them to be trivial.

\bigskip

\section{Homological conformal field theory} 

We briefly recall basics of homological conformal field theory. Let $\mathcal M_{g,p+q}$ be the moduli space of (not necessarily connected) Riemann surfaces of genus $g$ with a holomorphic map from a disjoint union of $p+q$ discs onto their disjoint images on the surface. Here the first $p$ discs are designated as incoming and the remaining $q$ discs as outgoing. There is a natural action of the product of symmetric groups $\Sigma_p\times \Sigma_q$ on $\mathcal M_{g,p+q}$ by relabeling incoming and outgoing holomorphic discs. Let $\mathcal C$ be a category such that the set of objects is the set of nonnegative integers $Ob(\mathcal C)=\mathbb N\cup\{0\}$, and for $p,q\in Ob(\mathcal C)$ the morphism set from $p$ to $q$ is $\mathcal C(p,q)=\coprod_{g\ge0}\mathcal M_{g,p+q}$. Disjoint union and sewing of Riemann surfaces give rise to operations:
\begin{align*}
\otimes&: \mathcal C(p,q)\times\mathcal C(p',q') \longrightarrow \mathcal C(p+p',q+q'),\\
\circ &: \mathcal C(q,r)\times \mathcal C(p,q) \longrightarrow \mathcal C(p,r). 
\end{align*}
With respect to the tensor law on objects given by $p\otimes q=p+q$, the category $\mathcal C$ has the structure of a strict symmetric monoidal category. A symmetric monoidal functor from the category $\mathcal C$ to the category of complex vector spaces is a conformal field theory \cite{Se}. Thus a conformal field theory is a representation theory of moduli spaces of Riemann surfaces. 

Let $\HC$ be a strict symmetric monoidal category whose set of objects is $\mathbb N\cup\{0\}$, and whose set of morphisms from $p$ to $q$ for $p,q\in Ob(\HC)$ is given by 
\begin{equation*}
\HC(p,q)=H_*(\mathcal C(p,q))=\bigoplus_{g\ge0} H_*(\mathcal M_{g,p+q}).
\end{equation*}
The composition of morphisms come from the gluing operation of Riemann surfaces 
\begin{equation*}
\circ: H_*\bigl(\mathcal C(q,r)\bigr)\otimes H_*\bigl(\mathcal C(p,q)\bigr) \longrightarrow H_*\bigl(\mathcal C(p,r)\bigr).
\end{equation*}
A homological conformal field theory (HCFT) is a symmetric monoidal functor $\mathcal F$ from the category $\HC$ to the category $\mathcal Gr_*$ of graded groups. For such a functor $\mathcal F$, the graded group $\mathcal F(1)=A_*$ comes equipped with linear maps 
\begin{equation*}
\mathcal F: H_*\bigl(\mathcal C(p,q)\bigr) \longrightarrow \text{Hom}(A_*^{\otimes p},A_*^{\otimes q})
\end{equation*}
for $p,q\ge0$, which are compatible with gluing of Riemann surfaces so that for homology classes $x\in H_*\bigl(\mathcal C(q,r)\bigr)$ and $y\in H_*\bigl(\mathcal C(p,q)\bigr)$, their composition $x\circ y\in H_*\bigl(\mathcal C(p,r)\bigr)$ satisfies the relation
\begin{equation*}
\mathcal F(x\circ y)=\mathcal F(x)\circ \mathcal F(y): A_*^{\otimes p} \longrightarrow A_*^{\otimes q} \longrightarrow  A_*^{\otimes r}.
\end{equation*}
Let $C$ be a Riemann sphere with one incoming and one outgoing disjoint holomorphic discs, and let $P$ be a Riemann sphere with two incoming and one outgoing holomorphic discs (a pair of pants). Their homology classes $c=[C]\in H_0(\mathcal M_{0,1+1})$ and $m=[P]\in H_0(\mathcal M_{0,2+1})$ are independent of conformal structures on $C$ and on $P$. The associated morphism $\mathcal F(c)=1 : A_* \rightarrow A_*$ is the identity on $A_*$, and $\mathcal F(m): A_*\otimes A_* \rightarrow A_*$ gives an associative product structure with unit on $A_*$.  

Let $\Hc$ be a strict symmetric monoidal category with the object set $\mathbb N\cup\{0\}$, and with the morphism set from $p$ to $q$ given by $\Hc(p,q)=H_0\bigl(\mathcal C(p,q)\bigr)$, which depends only on the topological type of surfaces. A functor $\mathcal F_0$ from the category $\Hc$ to the category $\mathcal Gr_*$ of graded groups is a topological quantum field theory (TQFT) \cite{At}. If $B_*=\mathcal F_0(1)$ is the graded group associated to the object $1\in Ob(\Hc)$, then $B_*$ has the structure of a Frobenius algebra. Note that every homological conformal field theory $\mathcal F$ restricts to a topological quantum field theory $\mathcal F_0$ by restricting the morphism set from $H_*\bigl(\mathcal C(p,q)\bigr)$ to $H_0\bigl(\mathcal C(p,q)\bigr)$. 

If the object set of the category $\HC$ and $\Hc$ is the set $\mathbb N$ of positive integers, then the corresponding homological conformal field theories and topological quantum field theories do not necessarily have units and counits, and these theories are called noncompact HCFT and noncompact TQFT, respectively. In the context of string topology $\mathbb H_*(LM)$ for finite dimensional manifold $M$, we only require $q$ to be positive. Such a theory is called a theory with positive boundary. In the context of string topology on the cohomology $H^*(LM)$ of simply connected infinite dimensional manifold $M$ whose based loop space $\Omega M$ has finite dimensional cohomology over a field $k$, we require both $p$ and $q$ to be positive. Thus we have a noncompact HCFT in this case. 

As mentioned in the introduction, the moduli space $\mathfrak M_{g,p+q}$ of connected genus $g$ Riemann surfaces with $p+q$ embedded holomorphic discs is homotopy equivalent to the classifying space $B\Gamma_{g,p+q}$ of the mapping class group $\Gamma_{g,p+q}$ when $p+q\ge1$. In this paper, we will mostly working with mapping class groups rather than moduli spaces of Riemann surfaces. Thus, we briefly describe some structures of the category $\HC$ including compositions of morphisms and actions of symmetric groups, in terms of surface diffeomorphisms and mapping class groups. 

For $i=1,2$, let $S_i$ be a smooth oriented (not necessarily connected) surface of genus $g_i$ with $p_i$ incoming and $q_i$ outgoing parametrized boundaries. Let $\text{Diff}^+(S_i,\partial)$ be the topological group of orientation preserving diffeomorphisms of $S_i$ fixing boundaries pointwise, and let $\Gamma(S_i)=\pi_0\bigl(\text{Diff}^+(S_i,\partial)\bigr)$ be the mapping class group of $S_i$. Suppose the number $q_1$ of outgoing boundaries of $S_1$ is the same as the number $p_2$ of incoming boundaries of $S_2$. Then we can sew two surfaces together to obtain a surface $S_2\#S_1$. Since self-diffeomorphisms of $S_1$ and $S_2$ can be combined together to a self-diffeomorphism of $S_2\#S_1$, we get a homomorphism $\circ: \Gamma(S_2)\times\Gamma(S_1) \rightarrow \Gamma(S_2\#S_1)$, which induces a map of their classifying spaces $\circ: B\Gamma(S_2)\times B\Gamma(S_1) \rightarrow B\Gamma(S_2\#S_1)$. The induced homology homomorphism is the composition of morphisms in the category $\HC$. 

Next we explain that $H_*(B\Gamma_{g,n})$ has a natural right $\Sigma_n$ action. In particular, $H_*(B\Gamma_{g,p+q})$ has a natural $\Sigma_p\times\Sigma_q$ action. To see this, let $F_{g,n}$ be a connected oriented smooth surface of genus $g$ with $n$ boundaries with parametrization given  by $\phi=(\phi_1,\phi_2,\dots,\phi_n):\coprod^n S^1 \xrightarrow{\cong}\partial F_{g,n}$. Let $D_{\partial}$ be the topological group of orientation preserving diffeomorphisms fixing boundaries pointwise, and let $D_{\Sigma_n}$ be the topological group of orientation preserving diffeomorphisms $f$ permuting parametrized boundaries in the sense that $f\circ \phi_i=\phi_{\tau(i)}$ for all $1\le i\le n$ and for some permutation $\tau\in \Sigma_n$. 
The group $D_{\partial}$ is a normal subgroup of $D_{\Sigma_n}$ with the quotient $\Sigma_n$. Let $\widetilde{\Gamma}_{g,n}=\pi_0(D_{\Sigma_n})$. Then we have exact sequences of groups: 
\begin{align*}
1 \longrightarrow D_{\partial} \longrightarrow  & D_{\Sigma_n} \longrightarrow \Sigma_n \longrightarrow 1, \\
1 \longrightarrow \Gamma_{g,n} \longrightarrow & \widetilde{\Gamma}_{g,n} 
\longrightarrow \Sigma_n \longrightarrow 1.
\end{align*}
Let $ED_{\Sigma_n} \rightarrow BD_{\Sigma_n}$ be the universal bundle where the group $D_{\Sigma_n}$ acts freely on $ED_{\Sigma_n}$ from the right. Since $n\ge1$, the natural projection $D_{\partial} \rightarrow \pi_0(D_{\partial})=\Gamma_{g,n}$ is a homotopy equivalence \cite{ES}, and hence we have a homotopy equivalence $BD_{\partial}=ED_{\Sigma_n}/D_{\partial}\simeq B\Gamma_{g,n}$. Since $\Sigma_n$ acts freely on $BD_{\partial}$ from the right, it acts on its homology $H_*(BD_{\partial})\cong H_*(B\Gamma_{g,n})$ from the right. 

There is another way to view $\Sigma_n$ action on the homology $H_*(B\Gamma_{g,n})=H_*(\Gamma_{g,n})$. This point of view is more relevant for the next section. For each $\tau\in\Sigma_n$, choose a diffeomorphism $f_{\tau}\in D_{\Sigma_n}$ whose restriction to boundaries gives the permutation $\tau$. Such $f_{\tau}$ is only unique up to $D_{\partial}$. The conjugation by $f_{\tau}$ from the right induces an automorphism of $D_{\partial}$, hence of its group of connected components $\Gamma_{g,n}$. However, this automorphism of $\Gamma_{g,n}$ can depend on the choice of $f_{\tau}$. Since inner automorphism of a group $\Gamma_{g,n}$ induces an identity on homology of $\Gamma_{g,n}$ (see for example \cite{B}, page 48), the conjugation action of $f_{\tau}$ on the homology $H_*(\Gamma_{g,n})$ depends only on $\tau\in\Sigma_n$. 

These two actions of $\Sigma_n$ on $H_*(B\Gamma_{g,n})$ are in fact the same. 
By regarding a $\mathbb Z[\widetilde{\Gamma}_{g,n}]$ free resolution of $\mathbb Z$ as a free resolution of $\mathbb Z[\Gamma_{g,n}]$ and considering its geometric realization, we can easily see that this conjugation action of $\Sigma_n$ on $H_*(\Gamma_{g,n})$ from the right coincides with the one induced by the $\Sigma_n$ free action above on the classifying space $BD_{\partial}$ from the right. 

\begin{remark}\label{trivial S_n action}
As it turns out that the $\Sigma_n$ action on $H_k(\Gamma_{g,n})$ is trivial for $g\ge 2k$ (See Lemma 3.3 in \cite{BT}). This is a consequence of Harer-Ivanov stability theorem. Since this fact is relevant in the next section, we explain the reason. For $g\ge 2k$ we have an onto map $H_k(\Gamma_{g,1}) \to H_k(\Gamma_{g,n})$ using Ivanov's reformulation of the stability theorem stated in Theorem \ref{Ivanov stability}. For any element $x\in H_k(\Gamma_{g,n})$, let $y$ be a cycle in the bar complex of $\text{Diff}^+(F_{g,1})$ representing $x$. The surface $F_{g,n}$ can be decomposed as $F_{g,n}=F_{g,1}\#F_{0,n+1}$. For any $\tau\in\Sigma_n$, let $f_{\tau}$ be a diffeomorphism of $F_{0,n+1}$ which induces the permutation on $n$ boundaries not used for sewing with $F_{g,1}$ and which is identity on the boundary used for sewing. Since diffeomorphisms appearing in the expression of the cycle $y$ and $f_{\tau}$ have disjoint support on $F_{g,n}$, conjugation action of $f_{\tau}$ on $y$ is trivial. Hence the action of $\Sigma_n$ on $H_k(\Gamma_{g,n})$ is trivial. This $\Sigma_n$-invariance in the stable range has an interesting consequence in terms of HCFT. Namely, in the stable range $g\ge 2k$, any element $x\in H_k(\Gamma_{g,p+q})$ defines an $\Sigma_p\times \Sigma_q$-invariant operation $\mathcal F(x): A_*^{\otimes p} \longrightarrow A_*^{\otimes q}$. Thus we have an operation between symmetric powers of $A_*$:
\begin{equation*}
\mathcal F(x):S^p(A_*) \longrightarrow S^q(A_*),
\end{equation*}
where the first $S^p(A_*)$ is the symmetric quotient of $A^{\otimes p}$, and the second $S^q(A_*)$ is the $\Sigma_q$-invariants in $A^{\otimes q}$. 
\end{remark}

Our final remark in this section is to point out that $\Sigma_p\times\Sigma_q$-equivariance of the closed string topology operation \eqref{string operation} and \eqref{cohomology string operation} essentially comes from the following strictly commutative diagram, where $(\sigma,\tau)\in\Sigma_p\times\Sigma_q$, and $F$ is a surface with $p+q$ parametrized boundaries. Here, the left and the middle horizontal maps are induced by restriction to incoming or outgoing boundaries of $F$. 
\begin{equation*}
\begin{CD}
BD_{\partial}\!\times \!(LM)^p @<<< ED_{\Sigma_n}\!\!\underset{D_{\partial}}\times\!\!\text{Map}(F,M) @>>> BD_{\partial}\!\times \!(LM)^q @>>> (LM)^q \\
@V{\cong}V{(\sigma, \tau)\times \sigma}V  @V{\cong}V{(\sigma,\tau)}V  @V{\cong}V{(\sigma,\tau)\times\tau}V  @V{\cong}V{\tau}V \\
BD_{\partial}\!\times \!(LM)^p @<<< ED_{\Sigma_n}\!\!\underset{D_{\partial}}\times\!\!\text{Map}(F,M) @>>> BD_{\partial}\!\times \!(LM)^q @>>> (LM)^q
\end{CD}
\end{equation*}
Note that $BD_{\partial}$ as well as the second space from the left admit free actions of the entire symmetric group $\Sigma_{p+q}$.

\bigskip

\section{Vanishing theorem for closed string topology (I)}\label{proof of theorem}

We carefully examine Harer's Stability Theorem from HCFT point of view. For this purpose, we fix a connected smooth oriented surface $F_{g,p+q}$ for each $g\ge0$ and for each $p,q$ with $p+q\ge1$. Let $T=F_{1,1+1}$ be a torus with one incoming and one outgoing parametrized boundaries. The surface resulting from sewing $T$ to the $i$-th incoming boundary of $F_{g,p+q}$ is denoted by $F_{g,p+q}\#_iT$ for $1\le i\le p$. Similarly, when we sew $T$ to the $j$-th outgoing boundary of $F_{g,p+q}$, the resulting surface is denoted by $T\#_jF_{g,p+q}$ for $1\le j\le q$. These are genus $g+1$ surfaces with $p+q$ boundaries, and there exist orientation preserving diffeomorphisms $h_i$ and $h_j$ from these surfaces to $F_{g+1,p+q}$: 
\begin{equation*}
F_{g,p+q}\#_iT \xrightarrow[\cong]{h_i} F_{g+1,p+q} \xleftarrow[\cong]{h_j} T\#_jF_{g,p+q}.
\end{equation*}
Both $h_i$ and $h_j$ are determined up to post-compositions with elements in $D_{g+1,\partial}=\text{Diff}^+(F_{g+1,p+q},\partial)$. 
Since both diffeomorphisms $f\in\text{Diff}^+(F_{g,p+q},\partial)$ and $g\in\text{Diff}^+(T,\partial)$ fix boundaries, they can be glued along a boundary to obtain a diffeomorphism $f\# g\in\text{Diff}^+(F_{g+1,p+q},\partial)$. By taking their isotopy classes, we obtain a homomorphism of groups and an associated homomorphism in homology:
\begin{align*}
\Phi_i &:\Gamma_{g,p+q}\times\Gamma_{1,1+1} \longrightarrow \Gamma_{g+1,p+q}, \\
(\Phi_i)_* &: H_*(\Gamma_{g,p+q})\otimes H_*(\Gamma_{1,1+1}) \longrightarrow  H_*(\Gamma_{g+1,p+q}), 
\end{align*}
given by $\Phi_i([f],[g])=[h_i\circ(f\#g)\circ h_i^{-1}]$. Since different choices of $h_i$ differ by elements of $D_{g+1,\partial}$, and since every inner automorphism induces identity on homology, the homology homomorphism $(\Phi_i)_*$ depends only on $i$. Similarly, if we glue the torus $T$ to $j$-th outgoing boundary of $F_{g,p+q}$, we obtain homomorphisms 
\begin{align*}
\Psi_j&: \Gamma_{1,1+1}\times \Gamma_{g,p+q} \longrightarrow \Gamma_{g+1,p+q},\\
(\Psi_j)_* &: H_*(\Gamma_{1,1+1})\otimes H_*(\Gamma_{g,p+q}) \longrightarrow H_*(\Gamma_{g+1,p+q}).
\end{align*}
Let $\varphi_i: \Gamma_{g,p+q} \rightarrow \Gamma_{g+1,p+q}$ be given by $\varphi_i(z)=\Phi_i(z,1)$ for $z\in \Gamma_{g,p+q}$, where $1\in\Gamma_{1,1+1}$ is the unit. Similarly we let $\psi_j: \Gamma_{g,p+q} \rightarrow \Gamma_{g+1,p+q}$ be defined by $\psi_j(z)=\Psi_j(1,z)$. The induced homology maps $(\varphi_i)_*$ and $(\psi_j)_*$ for $1\le i\le p$ and $1\le j\le q$ are Harer's stabilizing maps. These stabilizing maps depend on $i,j$, but only up to $\Sigma_{p+q}$-equivariance, as we show next. 

First, note that the mapping class group $\Gamma_{g,p+q}$ does not distinguish between incoming and outgoing boundaries. Thus any statement on homology of mapping class groups before applying HCFT functor must be independent of the distinction between incoming and outgoing boundaries. So for convenience, for $1\le j\le q$ let $\varphi_{p+j}=\psi_j$, and we write $T\#_jF_{g,p+q}$ as $F_{g,p+q}\#_jT$ for uniformity of notation. 

\begin{proposition}\label{transposition}
For $1\le i, j\le p+q$, the homomorphisms $(\varphi_i)_*$ and $(\varphi_j)_*$ are related by the right action of  the transposition $\tau_{ij}\in\Sigma_{p+q}$ as in the following diagram\textup{:}
\begin{equation}\label{transposition diagram}
\begin{CD}
H_k(B\Gamma_{g,p+q}) @>{(\varphi_i)_*}>> H_k(B\Gamma_{g+1,p+q}) \\
@V{\cong}V{\cdot\tau_{ij}}V    @V{\cong}V{\cdot\tau_{ij}}V  \\
H_k(B\Gamma_{g,p+q}) @>{(\varphi_j)_*}>> H_k(B\Gamma_{g+1,p+q}). 
\end{CD}
\end{equation}
In the stable range of $g\ge 2k$, the action of $\Sigma_{p+q}$ is trivial and all stabilizing homomorphisms are the same\textup{:} $(\varphi_i)_*=(\varphi_j)_*$ for $1\le i,j\le p+q$. 
\end{proposition}
\begin{proof} As before, we choose diffeomorphisms $h_i$ and $h_j$ to the surface $F_{g+1,p+q}$ as in the following diagram: 
\begin{equation*}
F_{g,p+q}\#_iT \xrightarrow[\cong]{h_i} F_{g+1,p+q} \xleftarrow[\cong]{h_j} F_{g,p+q}\#_jT.
\end{equation*}
Choose a diffeomorphism $u_{ij}$ of $F_{g,p+q}$ which switches the $i$-th and the $j$-th boundaries, and which fixes other boundaries pointwise. The map $u_{ij}$ is unique up to post and pre-composition with elements in $\text{Diff}^+(F_{g,p+q},\partial)$. The map $u_{ij}$ induces a diffeomorphism $u_{ij}\#1: F_{g,p+q}\#_iT \xrightarrow{\cong} F_{g,p+q}\#_jT$ which is identity on $T$ and switches the $i$-th and the $j$-th boundaries. Now two diffeomorphisms $h_j\circ(u_{ij}\#1)$ and $h_i$ differ by a self-diffeomorphism $v_{ij}$ of $F_{g+1,p+q}$ switching the $i$-th and the $j$-th boundaries. Thus we have $h_j\circ(u_{ij}\#1)=v_{ij}\circ h_i$. Since $h_i$ and $h_j$ are unique up to post-composition by elements in $\text{Diff}^+(F_{g+1,p+q},\partial)$, the map $v_{ij}$ is unique up to post and pre-composition with elements in $\text{Diff}^+(F_{g+1,p+q},\partial)$. Since $\varphi_i([f])=[h_i\circ(f\#1)\circ h_i^{-1}]$ for $[f]\in\Gamma_{g,p+q}$ for all $i$, it is straightforward to check the commutativity of the following diagram: 
\begin{equation*}
\begin{CD}
\Gamma_{g,p+q} @>{\varphi_i}>>  \Gamma_{g+1,p+q}  \\
@V{[u_{ij}]\circ(\ )\circ [u_{ij}]^{-1}}V{\cong}V
@V{\cong}V{[v_{ij}]\circ(\ )\circ [v_{ij}]^{-1}}V \\
\Gamma_{g,p+q} @>{\varphi_j}>>  \Gamma_{g+1,p+q}.
\end{CD}
\end{equation*}
As observed earlier, elements $[u_{ij}]$ and $[v_{ij}]$ are unique up to post and pre-composition with elements in $\Gamma_{g,p+q}$ and in $\Gamma_{g+1,p+q}$, respectively. Thus, on homology level, they induce unique maps, namely the action by the transposition $\tau_{ij}\in\Sigma_{p+q}$, and the homology commutative diagram \eqref{transposition diagram} follows. 

When we are in the stable range $g\ge 2k$, by Remark \ref{trivial S_n action} the action of the symmetric group $\Sigma_{p+q}$ is trivial. Thus, all the stabilizing maps $(\varphi_i)_*$ for $1\le i\le p+q$ are the same. This completes the proof. 
\end{proof} 

Let $\mathcal F:\HC \rightarrow \mathcal Gr_*$ be a HCFT with $\mathcal F(1)=A_*$, and let $\mathcal F_0:\Hc \rightarrow \mathcal Gr_*$ be the associated TQFT obtained by restriction from $\mathcal F$. For $x\in H_*(B\Gamma_{g,p+q})$, we compare HCFT operations $\mathcal F(x)$ and $\mathcal F\bigl((\varphi_i)_*(x)\bigr):A_*^{\otimes p} \rightarrow A_*^{\otimes q}$ associated to $x$ and $(\varphi_i)_*(x)$, where the latter belongs to $H_*(B\Gamma_{g+1,p+q})$. Let $T=F_{1,1+1}$ be as before, and let $t=[T]\in H_0(B\Gamma_{1,1+1})\cong\mathbb Z$ be the generator. We also consider similar questions for $\psi_j$'s for $1\le j\le q$. 

\begin{proposition}\label{stabilizing map}
For $1\le i\le p$ and $x\in H_k(B\Gamma_{g,p+q})$, we have 
\begin{equation*}
\mathcal F\bigl((\varphi_i)_*(x)\bigr)=\mathcal F(x)\circ(1\otimes\cdots\otimes 1\otimes \mathcal F_0(t)\otimes 1\otimes \cdots\otimes 1): A_*^{\otimes p} \longrightarrow A_*^{\otimes q}.
\end{equation*}
Here $\mathcal F_0(t): A_* \rightarrow A_*$ is the TQFT operator associated to the torus $T$, inserted at the $i$-th position. For $\psi_j$ with $1\le j\le q$, the corresponding formula is 
\begin{equation*}
\mathcal F\bigl((\psi_j)_*(x)\bigr)=(1\otimes\cdots\otimes 1\otimes \mathcal F_0(t)\otimes 1\otimes \cdots\otimes 1)\circ\mathcal F(x): A_*^{\otimes p} \longrightarrow A_*^{\otimes q},
\end{equation*}
where $\mathcal F_0(t)$ is inserted at the $j$-th position. 

In the stable range of $g\ge 2k$, the above two operations are the same and defines a map between symmetric powers\textup{:}
\begin{equation*}
\mathcal F\bigl((\varphi_i)_*(x)\bigr)=\mathcal F\bigl((\psi_j)_*(x)\bigr)
:S^p(A_*) \longrightarrow S^q(A_*),
\end{equation*}
for $1\le i\le p$ and $1\le j\le q$.
\end{proposition}
\begin{proof} The homology homomorphism $(\Phi_i)_*=\circ_i: H_*(B\Gamma_{g,p+q})\otimes H_*(B\Gamma_{1,1+1}) \rightarrow H_*(B\Gamma_{g+1,p+q})$ induced from a group homomorphism $\Phi_i$ is part of the composition of morphisms $\HC(p,q)\otimes \HC(p,p) \rightarrow \HC(p,q)$ in the category $\HC$. Thus, by gluing property of HCFT $\mathcal F$, for $x\in H_*(B\Gamma_{g,p+q})$ and $y\in H_*(B\Gamma_{1,1+1})$ we have 
\begin{equation*}
\mathcal F(x\circ_i y)=\mathcal F(x)\circ(1\otimes\cdots\otimes\mathcal F(y)\otimes\cdots\otimes 1): A_*^{\otimes p} \longrightarrow A_*^{\otimes q},
\end{equation*}
where $\mathcal F(y)$ is at the $i$-th position. Since $\varphi_i:\Gamma_{g,p+q} \rightarrow \Gamma_{g+1,p+q}$ is given by $\varphi_i(z)=\Phi_i(z,1)$, the induced map on classifying spaces is given by 
\begin{equation*}
B\varphi_i: B\Gamma_{g,p+q} \overset{\cong}\longrightarrow B\Gamma_{g,p+q}\times \{*\} 
\longrightarrow B\Gamma_{g.p+q} \times B\Gamma_{1,1+1} \xrightarrow{\Phi_i} B\Gamma_{g+1,p+q},
\end{equation*}
where $*\in B\Gamma_{1,1+1}$ is any point. Since $[*]=t=[T]\in H_0(B\Gamma_{1,1+1})$, for any element $x\in H_*(B\Gamma_{g,p+q})$, we have $(\varphi_i)_*(x)=(\Phi_i)_*(x\otimes t)=x\circ_i t$. Since $\mathcal F(t)=\mathcal F_0(t)$ by definition, we have the formula in the proposition for $\mathcal F\bigl((\varphi_i)_*(x)\bigr)$. 

The proof for $\psi_j$'s is similar. 

In the stable range, the action of $\Sigma_{p+q}$ is trivial by Remark \ref{trivial S_n action}, and we have operations on symmetric powers. By Proposition \ref{transposition}, elements $(\varphi_i)_*(x)$ and $(\psi_j)_*(x)$ for $x\in H_k(B\Gamma_{g,p+q})$ are the same in $H_k(B\Gamma_{g+1,p+q})$, and hence give rise to the same HCFT operation. This completes the proof. 
\end{proof} 
 
The commutativity diagram in Proposition \ref{transposition} implies that for $x\in H_k(B\Gamma_{g,p+q})$, we have $\mathcal F\bigl((\varphi_i)_*(x)\cdot\tau_{ij}\bigr)=\mathcal F\bigl((\varphi_j)_*(x\cdot\tau_{ij})\bigr)$. When $1\le i,j\le p$, this formula can be verified directly using Proposition \ref{stabilizing map} in view of $\Sigma_p$-equivariance of string operations as follows. For $a_1,a_2,\dots,a_p\in A_*$, we have
\begin{multline*}
\mathcal F\bigl((\varphi_i)_*(x)\cdot\tau_{ij}\bigr)(a_1\otimes a_2\otimes \cdots\otimes a_p)
=\mathcal F\bigl((\varphi_i)_*(x)\bigr)
\bigl(\tau_{ij}(a_1\otimes\cdots\otimes a_p)\bigr)\\
=(-1)^{\varepsilon}\mathcal F(x)(a_1\otimes\cdots\otimes \mathcal F_0(t)a_j\otimes\cdots\otimes a_i\otimes\cdots\otimes a_p),
\end{multline*}
where $\mathcal F_0(t)a_j$ is at the $i$-th position and $a_i$ is at the $j$-th position, and the sign $(-1)^{\varepsilon}$ is given by 
\begin{equation*}
\varepsilon=|a_i||a_j|+(|a_i|+|a_j|)(|a_{i+1}|+\cdots+|a_{j-1}|)+|\mathcal F_0(t)|(|a_1|+\cdots+|a_{i-1}|).
\end{equation*}
On the other hand, 
\begin{multline*}
\mathcal F\bigl((\varphi_j)_*(x\cdot\tau_{ij})\bigr) (a_1\otimes\cdots\otimes a_p)=\mathcal F(x\cdot\tau_{ij})(a_1\otimes\cdots \otimes a_i\otimes\cdots\otimes\mathcal F_0(t)a_j\otimes\cdots\otimes a_p)\\
=(-1)^{\varepsilon}\mathcal F(x)(a_1\otimes\cdots\otimes \mathcal F_0(t)a_j\otimes\cdots\otimes a_i\otimes\cdots\otimes a_p),
\end{multline*}
with the same $\varepsilon$ as above, where in the first line, $\mathcal F_0(t)a_j$ is at the $j$-th position, and in the second line $\mathcal F_0(t)a_j$ is at the $i$-th position. 

The corresponding formulas involving $\psi_j$'s can be checked  similarly using $\Sigma_q$-equivariance of the HCFT operation. In the mixed case, for $1\le i\le p$ and $1\le j\le q$, by Proposition \ref{transposition} we have $\mathcal F\bigl((\varphi_i)_*(x)\cdot\tau_{ij}\bigr)=\mathcal F\bigl((\psi_j)_*(x\cdot\tau_{ij})\bigr)$ for $\tau_{ij}\in\Sigma_{p+q}$. Since HCFT operations are only $\Sigma_p\times \Sigma_q$-equivariant, we cannot go any further in this case, although of course in the stable range we can eliminate $\tau_{ij}$ from this formula. 

If the HCFT $\mathcal F$ is defined for an object $p=0$, and the unit $1\in A_*$ with respect to the product structure $\mathcal F_0(m):
A_*^{\otimes 2}  \rightarrow A_*$ exists, then the homomorphism $\mathcal F_0(t):A_* \rightarrow A_*$ is simply given by multiplication by an element $\xi=\mathcal F_0(t)(1)\in A_*$. This is because capping one of the incoming boundaries of $F_{1,2+1}$ gives a surface $F_{1,1+1}=F_{1,0+1}\#F_{0,2+1}$. In the following commutative diagram in which vertical homomorphisms are induced by capping $p$ incoming boundaries by discs,
\begin{equation*}
\begin{CD}
H_*(\Gamma_{0,p+1}) @>>> H_*(\Gamma_{g,p+1}) \\
@VVV    @VVV \\
H_*(\Gamma_{0,1})=0  @>>> H_*(\Gamma_{g,1}),
\end{CD}
\end{equation*}
the right vertical homomorphism is an isomorphism for a sufficiently large genus $g$ by Harer's stability theorem. Since the top map factors through a trivial homology group due to $\Gamma_{0,1}=1$, it is a zero homomorphism. Using this observation, Tillmann \cite{Til} showed that if $A_*$ supports a HCFT with unit and $\xi=\mathcal F_0(t)(1)\in A_*$ is such that $\xi^n\not=0$ for all $n\ge1$, then the Batalin-Vilkovisky algebra structure in $A_*$ is trivial after localization $A_*[\xi^{-1}]$. Here note that if $\xi$ is a multiplicative torsion element with $\xi^m=0$ for some $m\ge1$, then $A_*[\xi^{-1}]=0$, and the situation is trivial. 

The situation we deal with is complementary to the above situation. Recall that $T$ denotes a Riemann surface with one incoming and one outgoing embedded discs, and $t=[T]\in H_0(B\Gamma_{1,1+1})\cong\mathbb Z$ is a generator. 

\begin{proposition} \label{trivial HCFT operations}
Let $\mathcal F:\HC \rightarrow \mathcal Gr_*$ be a homological conformal field theory with $\mathcal F(1)=A_*$, and let $\mathcal F_0=\mathcal F|_{H_0}$ be the associated topological quantum field theory. If $\mathcal F_0(t)=0: A_* \rightarrow A_*$, then operations $\mathcal F(x)$ associated to elements $x\in H_*(B\Gamma_{g,p+q})$ in the image of any stabilizing maps 
\begin{equation*}
(\varphi_i)_*, \ (\psi_k)_*: H_*(B\Gamma_{g-1,p+q}) \longrightarrow H_*(B\Gamma_{g,p+q}),\qquad 1\le i\le p, \quad 1\le k\le q,
\end{equation*}
are trivial. 
\end{proposition}
\begin{proof}
This is a direct consequence of the formula in Proposition \ref{stabilizing map}.
\end{proof} 

\begin{proof}[Proof of Vanishing Theorem \textup{(I)}] If $M$ is a finite dimensional smooth oriented closed manifold, then Cohen-Godin \cite{CG} and Godin \cite{Go2} showed that $\mathbb H_*(LM)$ supports a structure of HCFT with positive boundary ($q\ge1$). Previously we showed that all higher genus topological quantum field theory  operations in closed string topology are trivial \cite{T3}, \cite{T5}. In particular $\mathcal F_0(t)=0$ for the genus one case as an operator on $\mathbb H_*(LM)$. Hence by Proposition \ref{trivial HCFT operations}, all string operations associated to images of stabilizing maps are trivial. 

A proof for part (II) is given in the next section. 
\end{proof}

\bigskip

\section{Vanishing theorem for closed string topology (II)} \label{vanishing theorem (II)}

If $M$ is a simply connected infinite dimensional manifold with finite dimensional cohomology $H^*(\Omega M;k)$ for its based loop space with coefficients in a field $k$ of an arbitrary characteristic, then Chataur and Menichi \cite{CM} showed that $H^*(LM;k)$ carries the structure of a noncompact HCFT  requiring $p,q\ge1$. We will show that genus $1$ topological quantum field theory operator vanishes, $\mathcal F_0(t)=0$ in this noncompact HCFT. Also, we show that the Serre spectral sequence for the fibration $p:LM\to M$ collapses and the coproduct map in $H^*(LM;k)$ vanishes if $H^*(\Omega M;k)$ is an exterior algebra on odd degree generators, which is the case when $H^*(\Omega M;\mathbb Z)$ has no torsion elements of order divisible by the characteristic of the field $k$ (Corollary \ref{trivial coproduct}). Again, by Proposition \ref{trivial HCFT operations}, all string operations associated to elements in images of stabilizing maps are trivial. 

The main point here is that relevant transfer maps, the integration along the fiber, can be defined because the fiber $\Omega M$ of fibrations $f_{\text{out}}$, $g_{\text{out}}$, and $\overline{g}$ in \eqref{fibration diagram} and \eqref{fibration diagram 2} below behaves as a finite dimensional oriented manifold. In particular, applying the transfer map $f_{\text{out}}^!$ is essentially equivalent to taking Poincar\'e duality of the Pontrjagin product map in $\Omega M$. Thus if the cohomology of $LM$ can be written as a tensor product $H^*(LM;k)=H^*(M;k)\otimes H^*(\Omega M;k)$ (which is the case when $H^*(\Omega M;k)$ is an exterior algebra by Theorem \ref{cohomology of LM} below), then applying partial Poincar\'e duality along the cohomologically finite dimensional fiber gives $H^*(M;k)\otimes H_*(\Omega M;k)$ in which cohomology loop product is given by the cup product in $H^*(M;k)$ and the Pontrjagin product in $H^*(\Omega M;k)$. In the context of the previous section, transfer maps used for construction of string operations can be defined because the manifold $M$ itself is finite dimensional. If the homology of $LM$ can be written as a tensor product $H_*(LM)=H_*(M)\otimes H_*(\Omega M)$, then applying  partial Poincar\'e duality along the finite dimensional base $M$, we get $H^*(M)\otimes H_*(\Omega M)$, which is the loop homology algebra $\mathbb H_*(LM)$ for the finite dimensional $M$. Thus, in this sense, Chataur-Menichi construction of noncompact HCFT on the cohomology $H^*(LM;k)$ does produce none other than loop homology in case of an infinite dimensional simply connected space $M$.

First we briefly recall TQFT product $\mu$ and coproduct $\Phi$ in $H^*(LM;k)$ defined in \cite{CM}. Suppose the finite dimensional connected commutative Hopf algebra $H^*(\Omega M;k)$ is concentrated between the degrees $0\le *\le d$. In this case, by Hopf-Borel Theorem on the structure of Hopf algebras (\cite{MT}, Theorem 1.3 and Corollary 1.4 in Chapter 7), the finite dimensional connected commutative Hopf algebra $H^*(\Omega M;k)$ must be one of the following forms as an algebra, where $p$ is the characteristic of the field $k$:
\begin{enumerate}
\item[(i)] When $p=0$, $H^*(\Omega M;k)\cong \Lambda_k(x_1,x_2,\dots,x_{\ell})$, where $|x_i|$ is odd. 
\item[(ii)] When $p=2$, $H^*(\Omega M;k)\cong \bigotimes_{i=1}^rk[y_i]/(y_i^{2^{f_i}})$. 
\item[(iii)] When $p\not=0,2$, $H^*(\Omega M;k)\cong \Lambda_k(x_1,x_2,\dots,x_{\ell})\otimes \bigotimes_{j=1}^m\bigl(k[y_j]/(y_j^{p^{f_j}})\bigr)$, where $|x_i|$ is odd and $|y_j|$ is even. 
\end{enumerate}
In particular, $H^*(\Omega M;k)$ is a Poincar\'e duality algebra with an orientation class $[\Omega M]\in H^d(\Omega M;k)$ in the top degree. The cohomology loop product and loop coproduct maps 
\begin{align*}
\Phi&: H^*(LM;k) \longrightarrow H^*(LM;k)\otimes H^*(LM;k),\\
\mu&: H^*(LM;k)\otimes H^*(LM;k) \longrightarrow H^*(LM;k),
\end{align*}
are homomorphisms of degree $-d$ defined as follows. Let $F=F_{0,2+1}$ be a pair of pants with two incoming boundaries and one outgoing boundary.  Restriction to boundaries give rise to two fibrations:
\begin{equation}\label{two fibrations}
\begin{CD}
LM\times LM @<{g_{\text{out}}}<< \text{Map}\,(F,M) @>{g_{\text{in}}}>> LM.
\end{CD}
\end{equation}
Here we switched words ``in" and ``out" since in cohomology formulation, arrows are reversed as in \eqref{cohomology string operation}. Since the surface $F$ is homotopy equivalent to a graph \begin{tikzpicture}\draw (0,0) arc (0:180:0.1) -- ++(- 0.2,0) arc (0:360: 0.1) -- ++( 0.2,0) arc (180:360: 0.1);\fill (- 0.2,0) circle ( 0.035); \fill (- 0.4,0) circle ( 0.035);\end{tikzpicture} with an appropriate orientation, by replacing $\text{Map}(F,M)$ with a homotopy equivalent space $\text{Map}(\begin{tikzpicture}\draw (0,0) arc (0:180: 0.1) -- ++(- 0.2,0) arc (0:360: 0.1) -- ++( 0.2,0) arc (180:360: 0.1);\fill (- 0.2,0) circle ( 0.035); \fill (- 0.4,0) circle ( 0.035);\end{tikzpicture},M)$, we have the following commutative diagram where the square is a pull-back diagram of fibrations $g_{\text{out}}$ and $\overline{g}$ with fiber $\Omega M$:
\begin{equation}\label{fibration diagram}
\begin{CD}
LM\times LM @<{g_{\text{out}}}<< \text{Map}(\begin{tikzpicture}\draw (0,0) arc (0:180:0.1) -- ++(- 0.2,0) arc (0:360: 0.1) -- ++( 0.2,0) arc (180:360: 0.1);\fill (- 0.2,0) circle ( 0.035); \fill (- 0.4,0) circle ( 0.035);\end{tikzpicture},M) @>{g_{\text{in}}}>> LM \\
@V{p\times p}VV @V{q}VV @. \\
M\times M @<{\overline{g}}<< \text{Map}(I,M). @. 
\end{CD}
\end{equation}
The map $q$ above is the restriction to the interval between two circles, and the bottom map $\overline{g}$ is the evaluation map at end points of the unit interval $I=[0,1]$. The map $g_{\text{in}}$ is defined by an onto map $S^1 \to \begin{tikzpicture}\draw (0,0) arc (0:180: 0.1) -- ++(- 0.2,0) arc (0:360: 0.1) -- ++( 0.2,0) arc (180:360: 0.1);\fill (- 0.2,0) circle ( 0.035); \fill (- 0.4,0) circle ( 0.035);\end{tikzpicture}$ which maps the base point of $S^1$ to one of the vertices of the graph, and traces each circle of the graph once and traces the middle interval twice in opposite directions. See the description of the map $f_1$ in \eqref{correspondence maps} for details. Since $M\times M$ is simply connected by hypothesis, the map $\overline{g}: \text{Map}(I,M) \longrightarrow M\times M$ is an oriented fibration with fiber $\Omega M$. Namely, $\pi_1(M\times M)$ acts trivially on the orientation class $[\Omega M]\in H^d(\Omega M;k)$.  Consequently, the pull-back fibration $g_{\text{out}}: \text{Map}(\begin{tikzpicture}\draw (0,0) arc (0:180: 0.1) -- ++(- 0.2,0) arc (0:360: 0.1) -- ++( 0.2,0) arc (180:360: 0.1);\fill (- 0.2,0) circle ( 0.035); \fill (- 0.4,0) circle ( 0.035);\end{tikzpicture},M) \longrightarrow LM\times LM$ is also an oriented fibration with fiber $\Omega M$. 

Then we can consider the following transfer maps $g_{\text{out}}^!$ and $\overline{g}^!$ of degree $-d$, both integrations along the fiber: 
\begin{align*}
g_{\text{out}}^!&: H^*\bigl(\text{Map}(\begin{tikzpicture}\draw (0,0) arc (0:180: 0.1) -- ++(- 0.2,0) arc (0:360: 0.1) -- ++( 0.2,0) arc (180:360: 0.1);\fill (- 0.2,0) circle ( 0.035); \fill (- 0.4,0) circle ( 0.035);\end{tikzpicture},M);k\bigr) \longrightarrow H^*(LM;k)\otimes H^*(LM;k),\\
\overline{g}^!&: H^*\bigl(\text{Map}(I,M);k\bigr)=H^*(M;k) \longrightarrow H^*(M;k)\otimes H^*(M;k).  
\end{align*}
The coproduct map $\Phi$ in $H^*(LM;k)$ is defined in terms of the transfer map by  
\begin{equation*}
\Phi: H^*(LM;k) \xrightarrow{g_{\text{in}}^*} H^*\bigl(\text{Map}(\begin{tikzpicture}\draw (0,0) arc (0:180: 0.1) -- ++(- 0.2,0) arc (0:360: 0.1) -- ++( 0.2,0) arc (180:360: 0.1);\fill (- 0.2,0) circle ( 0.035); \fill (- 0.4,0) circle ( 0.035);\end{tikzpicture},M);k\bigr) \xrightarrow{g_{\text{out}}^!} H^*(LM;k)\otimes H^*(LM;k). 
\end{equation*}

To understand these transfer maps, we recall the Serre spectral sequence description of the integration along the fiber in a general form. Let $p:E \longrightarrow B$ be an oriented fibration with connected fiber $F$ such that the cohomology $H^*(F;k)$ with coefficient in a field $k$ is finite dimensional with a top degree orientation class $[F]\in H^d(F;k)\cong k$. Then the integration along the fiber $p^!$ is given by the following composition of maps and lowers cohomological degree by $d$: 
\begin{equation}\label{integration along fiber}
p^!: H^{n+d}(E;k)=D^{n,d} \xrightarrow{\text{onto}} E^{n,d}_{\infty}\subset E^{n,d}_2=H^n\bigl(B;H^d(F;k)\bigr) \xleftarrow[\cong]{[F]} H^n(B;k),
\end{equation}
where $H^{n+d}(E)=D^{0,n+d}\supset\cdots\supset D^{n,d}\supset D^{n+1,d-1}\supset\cdots\supset D^{n+d,0}\supset 0$ is a filtration on $H^{n+d}(E;k)$. Since $H^k(F;k)=0$ for $k>d$, we have $E_{\infty}^{n+d-k,k}=0$ for $k>d$. Thus, we have $H^{n+d}(E;k)=D^{0,n+d}=\cdots=D^{n,d}$, as in  \eqref{integration along fiber}.  

Here is a simple and useful criterion for vanishing of the transfer map $p^!$. This Lemma is used in Theorem \ref{cohomology of LM} below. 

\begin{lemma}\label{vanishing transfer} Let $p:E \longrightarrow B$ be an oriented fibration with connected fiber $F$ with an orientation class $[F]\in H^d(F;k)$ in the top degree. If $p^*:H^*(B;k) \longrightarrow H^*(E;k)$ is onto, then $p^!=0$. 
\end{lemma} 
\begin{proof} In terms of the Serre spectral sequence, the map $p^*$ is given by the following composition for an arbitrary $n$: 
\begin{equation*}
p^*: H^n(B;k)\cong H^n\bigl(B;H^0(F;k)\bigr)=E_2^{n,0} \xrightarrow{\text{onto}} E_{\infty}^{n,0}=D^{n,0}\subset D^{0,n}=H^n(E;k).
\end{equation*}
Thus, if $p^*$ is onto, then we have $D^{n,0}=D^{0,n}$, which is equivalent to $E_{\infty}^{*,q}=D^{*,q}/D^{*+1,q-1}=0$ for all $q\ge1$. In particular, $E_{\infty}^{*,d}=0$. Thus, $p^!=0$. This completes the proof. 
\end{proof}

Similarly, the intergration along the fiber in homology can be defined by the following composition of maps: 
\begin{equation}\label{homology transfer}
p_!:H_n(B;k)\cong H_n\bigl(B;H_d(F;k)\bigr)=E^2_{n,d} \xrightarrow{\text{onto}} E^{\infty}_{n,d}=D_{n,d}\subset D_{n+d,0}=H_{n+d}(E;k).
\end{equation}
We naturally expect that over a field coefficient $k$, the homology transfer $p_!$ and the cohomology transfer $p^!$ are dual to each other. Since we use this fact later, we quickly verify this. 

\begin{lemma}\label{dual transfers}  With the same hypothesis on the fibration $p:E\longrightarrow B$ as in Lemma \ref{vanishing transfer}, homology and cohomology transfer maps over a field $k$ are dual to each other, namely $(p_!)^*=p^!$. 
\end{lemma} 
\begin{proof} By comparing cohomology and homology transfers given in  \eqref{integration along fiber} and \eqref{homology transfer}, all we have to show is that the dual of homology $E^{\infty}$-terms are isomorphic to cohomology $E_{\infty}$-terms: $(E^{\infty}_{p,q})^*=E_{\infty}^{p,q}$. To see this, we have to recall the definition of homology and cohomology filtrations $\{D_{*,*}\}$ and $\{D^{*,*}\}$. These filtrations are defined in terms of an increasing subchain complexes $\{A_p\}$ of $C_*(E)$ by 
\begin{align*}
D_{p,q}&=\text{Im}\,[\iota_*:H_{p+q}(A_p;k) \longrightarrow H_{p+q}(E;k)]\\
D^{p,q}&=\text{Ker}\,[\iota^*:H^{p+q}(E;k) \longrightarrow H^{p+q}(A_{p-1};k)].
\end{align*}From this description, it is easy to see that homology and cohomology filtrations are related by 
\begin{equation}\label{cohomology filtration}
D^{p,q}\cong\bigl(H_{p+q}(E;k)/D_{p-1,q+1}\bigr)^*.
\end{equation}
By taking the dual of the following exact sequence,
\begin{equation*}
0 \longrightarrow E^{\infty}_{p,q} \longrightarrow H_{p+q}(E;k)/D_{p-1,q+1} \longrightarrow H_{p+q}(E;k)/D_{p,q} \longrightarrow 0,
\end{equation*}
and using \eqref{cohomology filtration}, we see that the above sequence becomes $0 \gets (E^{\infty}_{p,q})^*  \gets D^{p,q} \gets D^{p+1,q-1} \gets 0$. Hence $(E^{\infty}_{p,q})^*\cong E_{\infty}^{p,q}$. This completes the proof. 
\end{proof}

\begin{remark} Using the descriptions of $p^!$ and $p^*$ given above in terms of spectral sequences, it is easy to see that $p^!\circ p^*=0$ Similarly, we can show that $p_*\circ p_!=0$.
\end{remark}

Next, we describe cohomology loop product. By replacing a pair of pants $F_{0,1+2}$ with one incoming and two outgoing boundaries by a homotopy equivalent graph \begin{tikzpicture} \draw (0,0) arc (0:360: 0.1) -- ++(- 0.2,0); \fill (0,0) circle ( 0.035); \fill (- 0.2,0) circle ( 0.035);\end{tikzpicture} with an appropriate orientation, we can replace the diagram \eqref{two fibrations} with the following one, where the square is a pull-back diagram of fibrations $f_{\text{out}}$ and $\overline{g}$ with fiber $\Omega M$:
\begin{equation}\label{fibration diagram 2}
\begin{CD}
LM @<{f_{\text{out}}}<< \text{Map}(\begin{tikzpicture} \draw (0,0) arc (0:360: 0.1) -- ++(- 0.2,0); \fill (0,0) circle ( 0.035); \fill (- 0.2,0) circle ( 0.035);\end{tikzpicture},M) @>{f_{\text{in}}}>> LM\times LM \\
@V{p\times p_{\frac12}}VV @V{q}VV @. \\
M\times M @<{\overline{g}}<< \text{Map}(I,M). @.
\end{CD}
\end{equation}
Here $q$ is the restriction to the middle interval of the graph \begin{tikzpicture} \draw (0,0) arc (0:360: 0.1) -- ++(- 0.2,0); \fill (0,0) circle ( 0.035); \fill (- 0.2,0) circle ( 0.035);\end{tikzpicture}, $f_{\text{out}}$ is induced by the restriction to the outer circle, and $f_{\text{in}}$ is the restriction to boundaries of upper and lower half discs of the graph. See the description of maps $f_3$ and $f_4$  in \eqref{correspondence maps} for details.  
The cohomology loop product map $\mu$ in $H^*(LM;k)$ of degree $-d$ is then defined by
\begin{equation*}
\mu: H^*(LM;k)\otimes H^*(LM;k) \xrightarrow{f_{\text{in}}^*} H^*\bigl(\text{Map}(\begin{tikzpicture} \draw (0,0) arc (0:360: 0.1) -- ++(- 0.2,0); \fill (0,0) circle ( 0.035); \fill (- 0.2,0) circle ( 0.035);\end{tikzpicture};k\bigr) \xrightarrow{f_{\text{out}}^!} H^{*-d}(LM;k).
\end{equation*}
The product map $\mu$ is in general nontrivial, but the coproduct map $\Phi$ is often trivial. We will discuss two cases in which $\Phi$ is trivial. But before this we prove a general fact that the composition $\mu\circ\Phi$, the genus 1 TQFT operator, is always trivial over any coefficient field $k$. This is exactly what is needed for Vanishing Theorem in the introduction. 

\begin{theorem}\label{vanishing of genus one map} Let $M$ be simply connected with finite dimensional $H^*(\Omega M;k)$. Then the genus $1$ TQFT operator associated to $F_{1,1+1}$ is trivial. Namely, 
\begin{equation*}
\mu\circ\Phi=0: H^{*+2d}(LM;k) \xrightarrow{\Phi} H^{*+d}(LM;k)\otimes H^*(LM;k) \xrightarrow{\mu} H^*(LM;k).
\end{equation*}
\end{theorem}
\begin{proof}  We consider the following composition diagram of correspondences for the product $\mu$ and the coproduct $\Phi$, with renamed maps for convenience:

\begin{equation}\label{composition}
\xymatrix{
LM  &  \text{Map}(\begin{tikzpicture}\draw (0,0) arc (0:180: 0.1) -- ++(- 0.2,0) arc (0:360: 0.1) -- ++( 0.2,0) arc (180:360: 0.1);\fill (- 0.2,0) circle ( 0.035); \fill (- 0.4,0) circle ( 0.035); \end{tikzpicture}, M) \ar[l]_{f_1\ \ \ \ \ } \ar[r]^{\ \ f_2} & LM\times LM \\
& \text{Map}(\begin{tikzpicture} \draw (0,0) ellipse ( 0.1 and  0.035) ellipse ( 0.1 and  0.1); \fill ( 0.1,0) circle ( 0.035); \fill (- 0.1,0) circle ( 0.035);\end{tikzpicture}, M) \ar[u]_{f_3'} \ar[r]^{f_2'} 
\ar[ul]^{f_5=f_1\circ f_3'\ \ }  \ar[dr]_{f_6=f_4\circ f_2'\ \ } 
& \text{Map}(\begin{tikzpicture} \draw (0,0) arc (0:360: 0.1) -- ++(- 0.2,0); \fill (0,0) circle ( 0.035); \fill (- 0.2,0) circle ( 0.035);\end{tikzpicture}, M) \ar[u]_{f_3} \ar[d]^{f_4} \ar[u]_{f_3} \\
& & LM
}
\end{equation}
The coproduct map and the product map are given by $\Phi=f_2^!\circ f_1^*$ and $\mu=f_4^!\circ f_3^*$, respectively. 
We label elements in the mapping spaces $A\in \text{Map}(\begin{tikzpicture}\draw (0,0) arc (0:180: 0.1) -- ++(- 0.2,0) arc (0:360: 0.1) -- ++( 0.2,0) arc (180:360: 0.1);\fill (- 0.2,0) circle ( 0.035); \fill (- 0.4,0) circle ( 0.035); \end{tikzpicture} ,M)$, $B\in \text{Map}(\begin{tikzpicture} \draw (0,0) arc (0:360: 0.1) -- ++(- 0.2,0); \fill (0,0) circle ( 0.035); \fill (- 0.2,0) circle ( 0.035);\end{tikzpicture} ,M)$, and $C\in \text{Map}(\begin{tikzpicture} \draw (0,0) ellipse ( 0.1 and  0.035) ellipse ( 0.1 and  0.1); \fill ( 0.1,0) circle ( 0.035); \fill (- 0.1,0) circle ( 0.035);\end{tikzpicture},M)$ by labeling arcs and vertices of the above graphs using image arcs and image vertices in $M$ as follows, where $x,y$ are points in $M$ and $\alpha,\beta,\gamma,\dots$ denote arcs in $M$. For example, $A$ below represents a map $A:\begin{tikzpicture}\draw (0,0) arc (0:180: 0.1) -- ++(- 0.2,0) arc (0:360: 0.1) -- ++( 0.2,0) arc (180:360: 0.1);\fill (- 0.2,0) circle ( 0.035); \fill (- 0.4,0) circle ( 0.035); \end{tikzpicture} \to M$ such that two circles of the graph are mapped to loops $\rho$ and $\sigma$ in $M$, and the middle interval is mapped to an arc $\eta$ from a point $x$ to $y$ in $M$. Here, the point $x$ plays the role of the base point of the images in all three cases. 

\begin{center}
\begin{tikzpicture}[>=stealth]
\draw (0,0) arc (0:180:0.3) -- ++(-0.6,0) arc (0:360:0.3) -- ++(0.6,0) arc (180:360:0.3);
\fill (-1.2,0) circle ( 0.035); 
\fill (-0.6,0) circle ( 0.035); 
\path (-0.3,0.3) node[above] {$\sigma$};
\path (-1.5,0.3) node[above] {$\rho$};
\path (-0.9,0) node[above] {$\eta$};
\path (-0.6,0) node[right] {$y$};
\path (-1.2,0) node[left] {$x$};
\draw[->] (-0.85,0) -- ++(0.01,0);
\draw[->] (-0.35,0.3) -- ++(-0.01,0);
\draw[->] (-1.55,0.3) -- ++(-0.01,0);
\path (-1.9,0) node[left] {$A$:};

\path (2.5,0) coordinate (P1);
\draw (P1)++(0,0) arc (0:360:0.4) -- ++(-0.8,0);
\fill (P1)++(0,0) circle ( 0.035) ;
\fill (P1)++(-0.8,0) circle ( 0.035);
\path (P1)++(0,0) node[right] {$x$};
\path (P1)++(-0.8,0) node[left] {$y$};
\path (P1)++(-0.4,0.4) node[above] {$\alpha$};
\path (P1)++(-0.4,-0.53) node[above] {$\beta$};
\path (P1)++(-0.35,-0.08) node[above] {$\gamma$};
\draw [->] (P1)++(-0.45,0.4) -- ++(-0.01,0);
\draw [->] (P1)++(-0.35,-0.4) -- ++(0.01,0);
\draw[->] (P1)++(-0.35,0) -- ++(0.01,0);
\path (P1)++(-1.1,0) node[left] {$B$:};

\path (5,0) coordinate (P2);
\draw (P2)++(0,0) ellipse (0.5 and 0.3) 
            ellipse  (0.5 and 0.7);
\fill (P2)++(0.5,0) circle ( 0.035);
\fill (P2)++(-0.5,0) circle ( 0.035);
\path (P2)++(0,0.7) node[above] {$\alpha$};
\path (P2)++(0,0.25) node[above] {$\eta$};
\path (P2)++(0,-0.32) node[above] {$\gamma$};
\path (P2)++(0,-0.82) node[above] {$\beta$};
\path (P2)++(0.5,0) node[right] {$x$};
\path (P2)++(-0.5,0) node[left] {$y$};
\draw[->] (P2)++(-0.1,0.7) -- ++(-0.01,0);
\draw[->] (P2)++(-0.1,0.3) -- ++(-0.01,0);
\draw[->] (P2)++(0.1,-0.3) -- ++(0.01,0);
\draw[->] (P2)++(0.1,-0.7) -- ++(0.01,0);
\path (P2)++(-0.8,0) node[left] {$C$:};
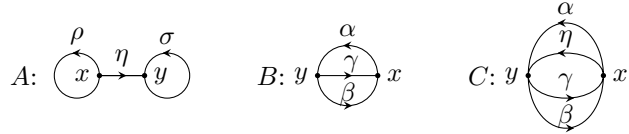
\path (2,-1) node[text width=11cm] {\textsc{Figure 1.} Descriptions of elements in three mapping spaces};
\end{tikzpicture}
\end{center}

In terms of this description, the top horizontal maps and the right vertical maps in the diagram \eqref{composition} are given in the following way: 
\begin{equation}\label{correspondence maps} 
f_1(A)=\rho\eta\sigma\eta^{-1}, \quad f_2(A)=(\rho, \sigma), \quad f_3(B)=(\alpha\gamma, \beta\gamma^{-1}), \quad f_4(B)=\alpha\beta. 
\end{equation}
The map $f_2'$ forgets the arc $\eta$, and the map $f_3'$ maps $C$ to $A$ with $\rho=\alpha\gamma$ and $\sigma=\beta\gamma^{-1}$. Finally, the diagonal maps in \eqref{composition} are given by 
\begin{equation*}
f_5(C)=\alpha\beta\eta\beta\gamma^{-1}\eta^{-1}, \qquad f_6(C)=\alpha\beta.
\end{equation*}
In the diagram \eqref{composition}, the maps $f_2, f_2', f_4$ are fibrations with fiber $\Omega M$ induced from $\overline{g}:\text{Map}(I,M) \longrightarrow M\times M$. Since the square in the diagram commutes, we have $f_3^*\circ (f_2)^!=(f_2')^!\circ(f_3')^*$. Hence 
\begin{equation*}
\mu\circ\Phi=(f_4^!\circ f_3^*)\circ (f_2^!\circ f_1^*)=(f_4\circ f_2')^!\circ(f_1\circ f_3')^*=f_6^!\circ f_5^*.
\end{equation*}
To understand the diagonal arrows of the diagram, consider the following pull-back diagram of fibrations: 
\begin{equation*}
\begin{CD}
LM\!\!\!\!\underset{M\times M}\times\!\!\!\! LM  @>{\pi_1}>> LM \\
@V{\pi_2}VV @V{(p_0, p_{\frac12})}VV  \\
LM @>{(p_0, p_{\frac12})}>> M\times M,
\end{CD}
\end{equation*}
where for a loop $\gamma:[0,1] \to M$ in $LM$, $p_0(\gamma)=\gamma(0)$ and $p_{\frac12}(\gamma)=\gamma(\frac12)$. 
The space $LM\times_{M\times M}LM$ consists of pairs of loops $(\ell_1,\ell_2)$ with the same base points and the same mid points, and maps $\pi_1,\pi_2$ are projections onto the first and the second components. In other words, any element in $LM\times_{M\times M}LM$ with base point $x$ and mid point $y$ is of the form $(\alpha\beta,\eta\gamma)$ where $\alpha,\eta$ are arcs from $x$ to $y$ and $\beta,\gamma$ are arcs from $y$ to $x$. Thus this space is exactly the same as the space $\text{Map}(\begin{tikzpicture} \draw (0,0) ellipse ( 0.1 and  0.035) ellipse ( 0.1 and  0.1); \fill ( 0.1,0) circle ( 0.035); \fill (- 0.1,0) circle ( 0.035);\end{tikzpicture},M)$. By homotopy equivalence, we can replace the space $LM\times_{M\times M}LM$ by a more familiar space. Let $LM\times_M LM\times_M LM$ be the space of triples of loops sharing the same base points. Consider maps  
\begin{align*}
h&: LM\!\!\!\!\underset{M\times M}\times\!\!\!\! LM \longrightarrow LM\underset{M}\times LM\underset{M}\times LM \\
\overline{h}&: LM\underset{M}\times LM\underset{M}\times LM \longrightarrow LM\!\!\!\!\underset{M\times M}\times\!\!\!\! LM,
\end{align*}
given by $h(\alpha\beta,\eta\gamma)=(\alpha\beta, \alpha\gamma, \eta\beta)$ and $\overline{h}(\alpha\beta,\xi_1,\xi_2)=\bigl(\alpha\beta,(\xi_2\beta^{-1})\cdot(\alpha^{-1}\xi_1)\bigr)$. See the diagram for  $C\in\text{Map}(\begin{tikzpicture} \draw (0,0) ellipse ( 0.1 and  0.035) ellipse ( 0.1 and  0.1); \fill ( 0.1,0) circle ( 0.035); \fill (- 0.1,0) circle ( 0.035);\end{tikzpicture},M)\cong LM\times_{M\times M}LM$ in Figure 1 above. We can easily check that these two maps are homotopy inverses to each other. Hence the diagonal part of the diagram \eqref{composition} can be replaced by the bottom line of the following diagram:  
\begin{equation*}
\begin{CD}
LM @<{f_5}<<  LM\!\!\!\!\underset{M\times M}\times\!\!\!\! LM @>{f_6}>> LM \\
@| @V{h}V{\simeq}V  @| \\
LM @<{f_7}<< LM\underset{M}\times LM\underset{M}\times LM @>{p_1}>>  LM,
\end{CD}
\end{equation*}
where $p_1(\delta,\xi_1,\xi_2)=\delta$, and $f_7(\delta,\xi_1,\xi_2)=\xi_1\xi_2\xi_1^{-1}\delta\xi_2^{-1}$. The right square of the above diagram strictly commutes and the left square commutes up to homotopy. Thus, 
\begin{equation*}
f_6^!\circ f_5^*=p_1^!\circ f_7^*: H^*(LM;k) \longrightarrow  H^*(LM\underset{M}\times LM\underset{M}\times LM;k) \longrightarrow  H^{*-2d}(LM;k).
\end{equation*}
To understand this map, we consider the dual homology map using Lemma \ref{dual transfers}: 
\begin{equation*}
(f_7)_*\circ(p_1)_!: H_*(LM;k) \longrightarrow H_{*+2d}(LM\underset{M}\times LM\underset{M}\times LM;k) \longrightarrow H_{*+2d}(LM;k).
\end{equation*}
We show that this composition is a zero map. To see this, let $w\in H_r(LM;k)$ be an arbitrary element, and let $W\subset M$ be a subspace of dimension at most $r$ such that an $r$-dimensional cycle representing $w$ is contained in $p^{-1}(W)\subset LM$ where $p:LM \rightarrow M$ is the base point projection. Let $L_WM=p^{-1}(W)$. Since $f_7$ preserves fibers over $M$, we have the following commutative diagram, where $\iota$'s are inclusion maps: 
\begin{equation*}
\begin{CD}
LM @<{f_7}<< LM\underset{M}\times LM\underset{M}\times LM @>{p_1}>> LM \\
@A{\iota}AA @A{\iota}AA  @A{\iota}AA  \\
L_WM @<{f_7^W}<< L_WM\underset{W}\times L_WM\underset{W}\times L_WM @>{p_1^W}>> L_WM 
\end{CD}
\end{equation*}
Let $w'\in H_r(L_WM;k)$ be an element such that $\iota_*(w')=w$. By the commutativity of the diagram, we have $(f_7)_*(p_1)_!(w)=\iota_*(f_7^W)_*(p_1^W)_!(w')$ where the element $(f_7^W)_*(p_1^W)_!(w')$ is in the group $H_{r+2d}(L_WM;k)$. 
Since $\dim W\le r$ and fiber $\Omega M$ has $k$-cohomological dimension $d$, we have $H_*(L_WM;k)=0$ in degrees $*>r+d$. Hence $H_{r+2d}(L_WM;k)=0$. Thus $(f_7)_*(p_1)_!(w)=0$. Since $w$ is arbitrary, it follows that the homology map $(f_7)_*(p_1)_!$ is a trivial map. Taking its dual, we see that the cohomology map $p_1^!f_7^*: H^{*+2d}(LM;k) \rightarrow H^*(LM;k)$ is also trivial. Consequently, we finally get  $\mu\circ\Phi=f_6^!f_5^*=p_1^!f_7^*=0$. 
\end{proof}

\begin{proof}[Proof of Vanishing Theorem \textup{(II)}] Theorem \ref{vanishing of genus one map} proves the vanishing of the TQFT operator $\mathcal F_0(t)=\mu\circ\Phi$ associated to the genus one surface $T_{\text{closed}}=F_{1,1+1}$ with one incoming and one outgoing boundaries. Again as in case (I),  Proposition \ref{trivial HCFT operations} proves the assertions in part (II) of Vanishing Theorem in the introduction. 
\end{proof} 

Next, we consider two cases in which the coproduct map $\Phi$ in the loop cohomology $H^*(LM;k)$ vanishes, which implies by dualizing that the product map in the loop homology $H_*(LM;k)$ vanishes, although the coproduct in $H_*(LM;k)$ is in general nontrivial. In the first case, we show that the existence of the unit with respect to the cohomology loop product implies the vanishing of the coproduct $\Phi$. Let $\iota:\Omega M \longrightarrow LM$ be the inclusion map of a fiber.

\begin{proposition} Let $M$ be simply connected with finite dimensional $H^*(\Omega M;k)$ concentrated in degrees $0\le *\le d$. Suppose there exists a unit $u\in H^d(LM;k)$ with respect to the cohomology loop product $\mu$. Then the followings hold\textup{:}

\textup{(i)} The coproduct map $\Phi$ is trivial. 

\textup{(ii)} The restriction of the unit to the fiber is the orientation class of the fiber. Namely, $\iota^*(u)=\{\Omega M\}\in H^d(\Omega M;k)$. 
\end{proposition}
\begin{proof} (i) In the diagram \eqref{composition}, the coproduct map $\Phi$ is given by $\Phi=f_2^!\circ f_1^*$. We consider the dual homology maps from a degree $0$ homology group: 
\begin{equation*}
(f_1)_*\circ (f_2)_!: H_0(LM\times LM;k) \longrightarrow H_d\bigl(\text{Map}(\begin{tikzpicture}\draw (0,0) arc (0:180: 0.1) -- ++(- 0.2,0) arc (0:360: 0.1) -- ++( 0.2,0) arc (180:360: 0.1);\fill (- 0.2,0) circle ( 0.035); \fill (- 0.4,0) circle ( 0.035); \end{tikzpicture}, M);k\bigr) \longrightarrow H_d(LM;k).
\end{equation*}
The generator of $H_0(LM\times LM;k)$ can be chosed to be $[(c_x,c_x)]$ for the constant loop $c_x$ at $x\in M$. Note that 
\begin{equation*}
f_1\circ f_2^{-1}\bigl((c_x,c_x)\bigr)=\{c_x\gamma c_x\gamma^{-1}\mid \gamma\in\Omega_xM\}\subset LM.
\end{equation*}
Let $\kappa:\Omega_xM \rightarrow LM$ be given by $\kappa(\gamma)=c_x\gamma c_x\gamma^{-1}$. Then $(f_1)_*(f_2)_!([(c_x,c_x)])=\kappa_*([\Omega_xM])$ in $H_d(LM;k)$, where $[\Omega_xM]$ is the orientation class in $H_d(\Omega_x M;k)$. Since obviously $\kappa$ is contractible, we have $\kappa_*([\Omega_xM])=0$. Thus, the composition of the above maps $(f_1)_*\circ (f_2)_!$ from degree $0$ homology is trivial. 

Let $z\in H^d(LM;k)$ be an arbitrary degree $d$ element. Since the dual of the coproduct map $\Phi$ is $(f_1)_*\circ (f_2)_!$, we have 
\begin{equation*}
\langle \Phi(z),[(c_x,c_x)]\rangle=\langle z,(f_1)_*(f_2)_!\bigl([(c_x,c_x)]\bigr)\rangle=\langle z,0\rangle=0. 
\end{equation*}
Hence it follows that $\Phi(z)=0\in H^0(LM\times LM;k)$. 

Now let $u\in H^d(LM;k)$ be the multiplicative unit in the loop cohomology $H^*(LM;k)$ so that for any $x\in H^*(LM;k)$, we have $u\cdot x=\mu(u\otimes x)=x$. Then the Frobenius relation implies that 
$\Phi(x)=\Phi(u\cdot x)=\Phi(u)\cdot x=0$, since $\Phi(u)=0$ because the degree of the unit $u$ is $d$. Hence the coproduct map $\Phi$ identically vanishes on the loop cohomology. 

(ii) We recall that the cohomology loop product map $\mu$ in $H^*(LM;k)$ is given by $\mu=(f_4)^!(f_3)^*$, where 
\begin{equation*}
\begin{CD}
LM\times LM @<{f_3}<< \text{Map}(\begin{tikzpicture} \draw (0,0) arc (0:360: 0.1) -- ++(- 0.2,0); \fill (0,0) circle ( 0.035); \fill (- 0.2,0) circle ( 0.035);\end{tikzpicture},M) @>{f_4}>> LM. 
\end{CD}
\end{equation*}
We examine the following dual homology map from a degree $0$ homology group:
\begin{equation*}
\begin{CD}
H_0(LM;k) @>{(f_4)_!}>> H_d\bigl(\text{Map}(\begin{tikzpicture} \draw (0,0) arc (0:360: 0.1) -- ++(- 0.2,0); \fill (0,0) circle ( 0.035); \fill (- 0.2,0) circle ( 0.035);\end{tikzpicture},M);k\bigr)  @>{(f_3)_*}>> H_d(LM\times LM;k).
\end{CD}
\end{equation*}
Since $M$ is simply connected, the free loop space $LM$ is connected. So a generator of $H_0(LM;k)$ can be chosen to be the class of the constant loop $[c_x]$ at $x\in M$. Then $(f_4)_!([c_x])=[\Omega_xM]\in H_d\bigl(\text{Map}(\begin{tikzpicture} \draw (0,0) arc (0:360: 0.1) -- ++(- 0.2,0); \fill (0,0) circle ( 0.035); \fill (- 0.2,0) circle ( 0.035);\end{tikzpicture},M);k \bigr)$ corresponding to the orientation class of the space of based loops obtained by mapping the middle interval of the graph \begin{tikzpicture} \draw (0,0) arc (0:360: 0.1) -- ++(- 0.2,0); \fill (0,0) circle ( 0.035); \fill (- 0.2,0) circle ( 0.035);\end{tikzpicture} into $M$, where the outer circle is mapped to a point $x\in M$ by $c_x$. Thus using the description of $f_3$ given in \eqref{correspondence maps}, the element $(f_3)_*(f_4)_!([c_x])$ is given by the homology class of the set $\{(\gamma,\gamma^{-1})\mid \gamma\in\Omega_xM\}\subset LM\times LM$. This is the image of the following composition: 
\begin{equation*}
\Omega_xM \xrightarrow{\phi} \Omega_x M\times \Omega_x M \xrightarrow{1\times S} \Omega_x M\times \Omega_x M \xrightarrow{\iota\times \iota} LM\times LM.
\end{equation*}
Thus $(f_3)_*(f_4)_!([c_x])=\iota_*[\Omega M]\otimes 1 + 1\otimes \iota_*S_*([\Omega M]) + \text{(other terms)}$. For $1\in H^0(LM;k)$, we have $\mu(u\otimes 1)=1$ because $u$ is the unit with respect to $\mu$. Since the dual of $\mu$ is $\mu^*=(f_3)_*(f_4)_!$ by Lemma \ref{dual transfers}, 
\begin{align*}
1&=\langle 1,[c_x]\rangle=\langle \mu(u\otimes 1),[c_x]\rangle=\langle u\otimes 1, (f_3)_*(f_4)_!([c_x])\rangle \\
&=\langle u\otimes 1, \iota_*[\Omega M]\otimes 1 + 1\otimes \iota_*S_*([\Omega M]) + \text{(other terms)}\rangle=\langle u, \iota_*([\Omega M])\rangle.
\end{align*}
Hence $\iota^*(u)([\Omega M])=1$. This completes the proof of (ii). 
\end{proof} 

In the second case in which $\Phi$ is trivial, we show that if $H^*(\Omega M;k)$ is an exterior algebra generated by finitely many odd degree elements, then the transfer map 
\begin{equation*}
f_2^!=g_{\text{out}}^!: H^*\bigl(\textup{Map}(\begin{tikzpicture}\draw (0,0) arc (0:180: 0.1) -- ++(- 0.2,0) arc (0:360: 0.1) -- ++( 0.2,0) arc (180:360: 0.1);\fill (- 0.2,0) circle ( 0.035); \fill (- 0.4,0) circle ( 0.035); \end{tikzpicture},M);k\bigr) \longrightarrow H^*(LM;k)\otimes H^*(LM;k). 
\end{equation*}
is already trivial. Consequently, the coproduct map $\Phi=f_2^!f_1^*$ is also trivial. To show this, we combine the Eilenberg-Moore spectral 	sequences and the Serre spectral sequences (see \cite{McC} for example). 

Let $p:E \to B$ be a fibration. In general, consider a diagram of pull-back fibrations and their induced cohomology maps: 
\begin{equation*}
\begin{CD}
E' @>{f}>> E \\
@V{p'}VV  @V{p}VV \\
B' @>{\overline{f}}>>  B
\end{CD}
\qquad \qquad
\begin{CD}
H^*(E';k) @<{f^*}<< H^*(E;k) \\
@A{(p')^*}AA  @A{p^*}AA \\
H^*(B';k)  @<{\overline{f}^*}<< H^*(B;k)
\end{CD}
\end{equation*}
If $p^*$ is onto, then under an appropriate condition on $H^*(F;k)$, Lemma \ref{vanishing transfer} implies $p^!=0$. What can we say about the other transfer $(p')^!$? We can use the Eilenberg-Moore spectral sequence, which is a second quadrant spectral sequence, to analyze the situation and to compute the cohomology $H^*(E';k)$. Its $E_2$-terms are given by 
\begin{equation}\label{Tor group}
E_2^{*,*}=\text{Tor}^{*,*}_{H^*(B)}\bigl(H^*(B'),H^*(E)\bigr).
\end{equation}
We consider a special case in which the Eilenberg-Moore spectral sequence collapses. 

\begin{lemma}\label{Eilenberg-Moore} Suppose $H^*(B';k)$ is a free $H^*(B;k)$-module. Then 
\begin{equation*}
H^*(E';k)\cong H^*(B';k)\!\!\!\!\underset{H^*(B)}\otimes\!\!\!\! H^*(E;k).
\end{equation*}
If furthermore, $p^*$ is onto and $I=\textup{Ker}\,p^*$, then 
\begin{equation*}
H^*(E';k)\cong H^*(B';k)/I',
\end{equation*}
where $I'=I\cdot H^*(B';k)$ is the extension of the ideal $I$ to an ideal in $H^*(B';k)$. 
\end{lemma} 
\begin{proof} Since $H^*(B';k)$ is a free module over $H^*(B;k)$, $\text{Tor}^{-p,q}=0$ for $p>0$ in \eqref{Tor group}. Thus the Eilenberg-Moore spectral sequence collapses and 
\begin{equation*}
H^*(E')\cong\text{Tor}^{0,*}_{H^*(B)}\bigl(H^*(B'),H^*(E)\bigr)=H^*(B')\!\!\!\!\underset{H^*(B)}\otimes\!\!\!\! H^*(E).
\end{equation*}
The second part follows from this. 
\end{proof} 

The pull-back diagram of fibrations relevant to us is the diagram \eqref{fibration diagram} which we reproduce here for convenience. 
\begin{equation*}
\begin{CD}
\text{Map}(\begin{tikzpicture}\draw (0,0) arc (0:180: 0.1) -- ++(- 0.2,0) arc (0:360: 0.1) -- ++( 0.2,0) arc (180:360: 0.1);\fill (- 0.2,0) circle ( 0.035); \fill (- 0.4,0) circle ( 0.035); \end{tikzpicture},M) @>{q}>> \text{Map}(I,M)\simeq M \\
@VV{f_2=g_{\text{out}}}V @V{\overline{g}=p_0\times p_1}V{\simeq \phi}V \\
LM\times LM @>{p\times p}>> M\times M
\end{CD}
\end{equation*}
Here, $\phi:M \to M\times M$ is the diagonal map. 
Since fibrations $q$ and $p\times p$ have sections, induced cohomology maps $g^*$ and $(p\times p)^*$ are injective. Since the induced cohomology map $\overline{g}^*=\phi^*:H^*(M;k)\otimes H^*(M;k) \longrightarrow H^*(M;k)$ is nothing but the cup product map, it is an onto map. Let $J=\text{Ker}\,\overline{g}^*$ be the kernel of the cup product map. 
If $H^*(LM;k)$ is a free $H^*(M;k)$-module, then Lemma \ref{Eilenberg-Moore} implies 
\begin{equation}\label{cohomology of mapping space}
H^*\bigl(\text{Map}(\begin{tikzpicture}\draw (0,0) arc (0:180: 0.1) -- ++(- 0.2,0) arc (0:360: 0.1) -- ++( 0.2,0) arc (180:360: 0.1);\fill (- 0.2,0) circle ( 0.035); \fill (- 0.4,0) circle ( 0.035); \end{tikzpicture},M);k\bigr)\cong H^*(LM\times LM;k)/J'
\end{equation}
and $g_{\text{out}}^*$ is onto, where $J'$ is an extension of $J$. By Lemma \ref{vanishing transfer}, we see that the transfer map $g_{\text{out}}^!$ is trivial. This proves the second part of the next theorem. 

\begin{theorem}\label{cohomology of LM} Let $k$ be a field of any characteristic. Let $M$ be simply connected such that 
\begin{equation}\label{cohomology of Omega M} 
H^*(\Omega M;k)\cong\Lambda_k(x_1,x_2,\dots,x_{\ell}),
\end{equation}
an exterior algebra on odd degree generators.  Then the following statements hold. 

\textup{(1)} The cohomology of $M$ is given by $H^*(M;k)\cong k[y_1,y_2,\dots,y_{\ell}]$, a polynomial algebra on even degree generators with $|y_i|=|x_i|+1$ for $1\le i\le \ell$, and the cohomology of $LM$ and $\textup{Map}(\begin{tikzpicture}\draw (0,0) arc (0:180: 0.1) -- ++(- 0.2,0) arc (0:360: 0.1) -- ++( 0.2,0) arc (180:360: 0.1);\fill (- 0.2,0) circle ( 0.035); \fill (- 0.4,0) circle ( 0.035); \end{tikzpicture},M)$ are given by 
\begin{align*}
H^*(LM;k)&\cong H^*(M;k)\otimes H^*(\Omega M;k) \\
    &\cong k[y_1,y_2,\dots,y_{\ell}]\otimes \Lambda_k(x_1,x_2,\dots,x_{\ell}), \\
H^*\bigl(\textup{Map}(\begin{tikzpicture}\draw (0,0) arc (0:180: 0.1) -- ++(- 0.2,0) arc (0:360: 0.1) -- ++( 0.2,0) arc (180:360: 0.1);\fill (- 0.2,0) circle ( 0.035); \fill (- 0.4,0) circle ( 0.035); \end{tikzpicture},M);k\bigr) &\cong H^*(M;k)\otimes H^*(\Omega M;k)\otimes H^*(\Omega M;k).
\end{align*}

\textup{(2)} The transfer map for the fibration $g_{\textup{out}}$ vanishes. Namely, 
\begin{equation*}
g_{\text{out}}^!=0: H^*\bigl(\textup{Map}(\begin{tikzpicture}\draw (0,0) arc (0:180: 0.1) -- ++(- 0.2,0) arc (0:360: 0.1) -- ++( 0.2,0) arc (180:360: 0.1);\fill (- 0.2,0) circle ( 0.035); \fill (- 0.4,0) circle ( 0.035); \end{tikzpicture},M);k\bigr) \longrightarrow H^*(LM;k)\otimes H^*(LM;k). 
\end{equation*}
Consequently, the coproduct map $\Phi=g_{\textup{out}}^!g_{\textup{in}}^*$ in $H^*(LM;k)$ also vanishes. Here the transfer $g_{\textup{out}}^!$ can be identified with 
\begin{multline*}
g_{\textup{out}}^!=\overline{g}^!\otimes 1\otimes 1: H^*(M;k)\otimes H^*(\Omega M;k)\otimes H^*(\Omega M;k) \\
\longrightarrow H^*(M;k)\otimes H^*(M;k)\otimes H^*(\Omega M;k)\otimes H^*(\Omega M;k). 
\end{multline*}
\end{theorem}

Over the rationals $k=\mathbb Q$, if $H^*(\Omega M;\mathbb Q)$ is finite dimensional, then it must be an exterior algebra on finitely many odd degree generators by the classical Hopf theorem. Thus hypothesis \eqref{cohomology of Omega M} is a natural one. See also Remark \ref{torsion elements and exterior algebra}. In this case, the rational cohomology $H^*(LM;\mathbb Q)$ can also be computed using minimal models. We emphasize that in Theorem \ref{cohomology of LM}, the characteristic of the field $k$ is arbitrary. 

The cohomology $H^*(LM;k)$ in part (1) of Theorem \ref{cohomology of LM} can be directly computed using Eilenberg-Moore spectral sequence. See for example \cite{Sm} where characteristic zero case is discussed. General characteristic $p$ case can be dealt with in the similar way. However, we found it more interesting and instructive (at least to the author) to show that the Serre spectral sequence for the fibration $\Omega M \to LM \xrightarrow{p} M$ collapses under our hypothesis on $H^*(\Omega M;k)$ given in \eqref{cohomology of Omega M}. See Remark \ref{differentials} at the end of this section for a motivation.  

We first discuss the structure of Serre spectral sequences of related fibrations. Let $x_0\in M$ be a base point, and let $\Omega M \to PM \xrightarrow{p} M$ be the path fibration, where $PM=\{\gamma:[0,1] \to M\mid \gamma(0)=x_0\}$ is the path space starting at $x_0$, and $p(\gamma)=\gamma(1)$. When the cohomology of the based loop space is given as in \eqref{cohomology of Omega M}, the Borel transgression theorem (\cite{MT}, Chapter 7, Theorem 2.9) tells us that the exterior algebra generators can be chosen to be transgressive so that 
\begin{equation*}
H^*(M;k)\cong k[y_1,y_2,\dots,y_{\ell}],
\end{equation*}
where $y_i$ is the image of $x_i$ under the transgression so that $|y_i|=|x_i|+1$ for $1\le i\le \ell$. In algebras $H^*(M;k)$ and $H^*(\Omega M;k)$, let $I(r)$ be the ideal generated by generators of degree less than $r$. We order generators $x_i$'s so that $|x_1|\le|x_2|\le\cdots\le|x_{\ell}|$, and let $|x_i|=r_i-1$ with $r_i=|y_i|$ even for $1\le i\le \ell$. The Serre spectral sequence $\{E_r^{*,*}\}$ for the path fibration $p:PM \to M$ is of the form 
\begin{equation}\label{spectral sequence for path fibration}
E_r^{*,*}=H^*(M;k)/I(r)\otimes H^*(\Omega M;k)/I(r-1),
\end{equation}
and all nonzero differentials in the spectral sequence are consequences of $d_{r_i}(x_i)=y_i$ for $1\le i\le \ell$, where $d_{r_i}: E_{r_i}^{0,r_i-1} \longrightarrow E_{r_i}^{r_i,0}$. 

Next we examine the Serre spectral sequence for the diagonal map  $\phi: M\longrightarrow M\times M$ regarded as a fibration. Let $J$ be the kernel of the cup product map $\phi^*:H^*(M)\otimes H^*(M) \longrightarrow  H^*(M)$. It can be easily checked that 
\begin{equation}\label{ideal J}
J=(y_1\otimes 1-1\otimes y_1, \ y_2\otimes 1-1\otimes y_2,\ \dots,\  y_{\ell}\otimes 1-1\otimes y_{\ell})\subset H^*(M\times M;k).
\end{equation}
Let $J(r)$ be the subideal of $J$ generated by generators of $J$ given in \eqref{ideal J} of degree less than $r$. 

\begin{lemma}\label{spectral sequence for g-bar} Let $M$ be simply connected with $H^*(\Omega M;k)$ given as in \eqref{cohomology of Omega M}. Let $\{E_r^{*,*}\}$ be the Serre spectral sequence for the fibration $\textup{Map}(I,M) \xrightarrow{\overline{g}=p_0\times p_1} M\times M$ with fiber $\Omega M$. Then its $E_r$-page is given by 
\begin{equation}\label{E_r}
E_r^{*,*}\cong H^*(M\times M;k)/J(r)\otimes H^*(\Omega M;k)/I(r-1),
\end{equation}
and all differentials are consequences of 
\begin{equation*}
d_{r_i}(x_i)=1\otimes y_i-y_i\otimes 1\in H^*(M\times M;k)/J(r_i), 
\end{equation*}
where $d_{r_i}:E_{r_i}^{0,r_i-1} \longrightarrow E_{r_i}^{r_i,0}$ for $1\le i\le \ell$. 
\end{lemma} 
\begin{proof} To examine the spectral sequence $\{E_r^{*,*}\}$, we compare it with related spectral sequences for path fibrations $PM \to M$ and $\overline{P}M \to M$, where $\overline{P}M=\{\gamma:[0,1] \to M\mid \gamma(1)=x_0\}$ is the space of paths ending at $x_0$. We have the following diagram 
\begin{equation*}
\begin{CD}
PM @>{h_1}>> \text{Map}(I,M) @<{h_2}<< \overline{P}M \\
@V{p_1}VV  @V{\overline{g}=p_0\times p_1}VV  @V{p_0}VV \\
\{x_0\}\times M @>{\bar{h}_1}>> M\times M @<{\bar{h}_2}<< M\times \{x_0\}, 
\end{CD}
\end{equation*}
where $h_1$ and $h_2$ are inclusions.
Let $S:PM \xrightarrow{\cong} \overline{P}M$ be the inverse path homeomorphism given by $S(\gamma)(t)=\gamma(1-t)$ for $t\in[0,1]$. Since the composite $\Omega M \xrightarrow{\phi} \Omega M\times \Omega M \xrightarrow{1\times S} \Omega M\times \Omega M \xrightarrow{m} \Omega M$, where $m$ is the loop multiplication map, is homotopic to the constant map, for any primitive element $z\in H^*(\Omega M;k)$, we have $S^*(z)=-z$. Since $H^*(\Omega M;k)$ is an exterior algebra generated by odd degree elements, by Hopf-Samelson Theorem (\cite{MT}, Chapter 7, corollary 1.13) $H^*(\Omega M;k)$ is primitively generated. Hence $S^*$ acts as $-1$ on odd degree elements. In particular, $S^*(x_i)=-x_i$ for $1\le i\le \ell$. (This argument is needed since in \eqref{cohomology of Omega M}, we did not assume that $x_i$'s are primitive.) Let $\{{}'E_r^{*,*}\}$ and $\{{}''E_r^{*,*}\}$ be Serre spectral sequences for the fibrations $PM\to M$ and $\overline{P}M \to M$. We discussed the spectral sequence $\{{}'E_r^{*,*}\}$ earlier in \eqref{spectral sequence for path fibration}. Using the isomorphism of spectral sequences $S^*:{}''E_r^{*,*} \xrightarrow{\cong}  {}'E_r^{*,*}$, we see that the differential $d_{r_i}'(x_i)=y_i$ in ${}'E_r^{*,*}$ translates to a differential $d_{r_i}''(x_i)=d_{r_i}'\bigl(S^*(x_i)\bigr)=-y_i$ in $\{{}''E_r^{*,*}\}$ for $1\le i\le \ell$, and all differentials in $\{{}''E_r^{*,*}\}$ are consequences of the above differentials. 

We prove Lemma \ref{spectral sequence for g-bar} by induction on $r\ge2$. When $r=2$, we have 
\begin{equation*}
E_2^{*,*}=H^*\bigl(M\times M;H^*(\Omega M;k)\bigr)=H^*(M\times M;k)\otimes H^*(\Omega M;k),
\end{equation*}
since the local system is trivial because $M$ is simply connected. Assume that $E_r$-page is given by \eqref{E_r}. Suppose $H^*(\Omega M;k)/I(r-1)\cong\Lambda_k(x_i,x_{i+1},\dots,x_{\ell})$ with $|x_i|\ge r-1>|x_{i-1}|$. We consider two cases. If $|x_i|>r-1$, then two consecutive degrees in $H^*(\Omega M)/I(r-1)$ with nontrivial groups are at least $|x_i|>r-1$ apart. Hence the differential $d_r:E_r^{*,*} \longrightarrow E_r^{*+r,*-r+1}$ must be trivial, and we have $E_{r+1}^{*,*}=E_r^{*,*}=H^*(M\times M)/J(r+1)\otimes H^*(\Omega M)/I(r)$ since $J(r+1)=J(r)$ and $I(r)=I(r-1)$ in this case. 

Suppose $|x_i|=r-1$ or $r=r_i$. We consider maps between spectral sequences induced by $h_1$ and $h_2$: $h_1^*: E_r^{*,*} \rightarrow {}'E_r^{*,*}$ and $h_2^*: E_r^{*,*} \rightarrow {}''E_r^{*,*}$.
We recall that ${}'E_r^{*,*}$ and ${}''E_r^{*,*}$ are given by 
\begin{equation*}
{}'E_r^{*,*}\cong H^*(M)/I(r)\otimes H^*(\Omega M)/I(r-1)\cong {}''E_r^{*,*}.
\end{equation*}
See \eqref{spectral sequence for path fibration}. Applying $h_1^*$ to the differential $d_r(x_i)$, we get $\bar{h}_1^*\bigl(d_r(x_i)\bigr)=d_r'(x_i)=y_i$. Similarly, applying $h_2^*$, we get $\bar{h}_2^*\bigl(d_r(x_i)\bigr)=d_r''(x_i)=-y_i$. Hence $d_r(x_i)$ is of the form $d_r(x_i)=1\otimes y_i-y_i\otimes 1+(\text{decomposables})$ in $H^*(M\times M)/J(r_i)$. Since the spectral sequence $\{E_r^{*,*}\}$ converges to $H^*\bigl(\text{Map}(I,M)\bigr)\cong H^*(M\times M)/J$, the decomposable elements in the above must lie in $J(r_i)$. Thus 
\begin{equation*}
d_{r_i}(x_i)=1\otimes y_i-y_i\otimes 1 \in H^*(M\times M)/J(r_i).
\end{equation*}
If there are other generators of $H^*(\Omega M)/I(r_i-1)$ of degree $r_i-1$, we can apply the same argument and we get corresponding results for their differentials. Hence 
\begin{equation*}
E_{r_i+1}^{*,*}=H^*(M\times M)/J(r_i+1)\otimes H^*(\Omega M)/I(r_i).
\end{equation*}
This completes the inductive step and the proof is complete. 
\end{proof} 

\begin{proof}[Proof of Theorem \ref{cohomology of LM}] We show that the Serre spectral sequence $\{E_r^{*,*}\}$ for the fibration $\Omega \to LM \xrightarrow{p} M$ collapses at $E_2$-term. Thus, suppose $E_r^{*,*}=E_2^{*,*}$ for some $r\ge2$. We show that $d_r=0$ on $E_r^{*,*}$, implying that $E_{r+1}^{*,*}=E_r^{*,*}$. Since $E_r^{*,*}=E_2^{*,*}=E_2^{*,0}\otimes E_2^{0,*}$, by derivation property of the differential, we only have to show that $d_r=0$ on $E_r^{0,*}=H^*(\Omega M;k)$. To see this, we consider the following pull-back diagram of fibrations:
\begin{equation}\label{pull-back for LM}
\begin{CD}
LM @>{h}>> \text{Map}(I,M) \\
@V{p}VV @V{\overline{g}=p_0\times p_1}VV \\
M @>{\phi}>> M\times M.
\end{CD}
\end{equation}
Let $\{{}'E_r^{*,*}\}$ be the Serre spectral sequence for the fibration $\overline{g}:\text{Map}(I,M) \longrightarrow M\times M$ described in Lemma \ref{spectral sequence for g-bar}. The map $h$ induces a map of spectral sequences $h^*: {}'E_r^{*,*} \longrightarrow E_r^{*,*}$. 
Suppose ${}'E_r^{0,*}=H^*(\Omega M;k)/I(r-1)=\Lambda_k(x_i,x_{i+1},\dots,x_{\ell})$, where $|x_i|\ge r-1>|x_{i-1}|$. 

We consider two cases. Suppose $|x_i|>r-1$. Since $d_r'$ is trivial on $x_i, x_{i+1}, \dots, x_{\ell}$ in the spectral sequence ${}'E_r^{0,*}$ in view of Lemma \ref{spectral sequence for g-bar}, mapping this relation to $E_r^{0,*}$ via $h^*$, we get $d_r(x_j)=0$ for $i\le j\le\ell$. By degree reason, $d_r$ is trivial on $x_1,x_2,\dots, x_{i-1}$, where $|x_1|\le\cdots\le|x_{i-1}|<r-1$. Hence $d_r=0$ on $E_r^{0,*}$. As remarked above, the derivation property of $d_r$ implies that $d_r=0$ on the entire $E_r^{*,*}$. 

Next, suppose $|x_i|=\cdots=|x_{k}|=r-1$ and $|x_{k+1}|>r-1$ for some $i\le k\le \ell$. In this case, $r=r_i$. Since $d_{r_i}'(x_j)=1\otimes y_j-y_j\otimes 1\in {}'E_{r_i}^{r_i,0}$ for $i\le j\le k$, mapping this relation by $h^*$, we get $d_{r_i}(x_j)=0\in E_{r_i}^{r_i,0}$ for $i\le j\le k$, since $h^*$ on the base spaces is simply the cup product $\phi^*$ (see diagram  \eqref{pull-back for LM}). By Lemma \ref{spectral sequence for g-bar}, for $x_j$ with $j>k$, we have $d_{r_i}'(x_j)=0$, which in turn implies that $d_{r_i}(x_j)=0$ for $k< j\le \ell$. By degree reason, $d_{r_i}$ is trivial on $x_1,x_2,\dots,x_{i-1}$. Hence again $d_{r_i}$ is trivial on $E_{r_i}^{0,*}$. Hence derivation property of differential implies that $d_{r_i}=0$ on the entire $E_{r_i}^{*,*}$. This completes the inductive step and we have proved the formula for the cohomology of $LM$. 

The cohomology of $\textup{Map}(\begin{tikzpicture}\draw (0,0) arc (0:180: 0.1) -- ++(- 0.2,0) arc (0:360: 0.1) -- ++( 0.2,0) arc (180:360: 0.1);\fill (- 0.2,0) circle ( 0.035); \fill (- 0.4,0) circle ( 0.035); \end{tikzpicture},M)$ then follows from \eqref{cohomology of mapping space}. This completes the proof of part (1). 

(2) Although the proof of part (2) was discussed right before the statement of Theorem \ref{cohomology of LM}, we give an alternate proof from a different point of view. 

The square diagram in \eqref{fibration diagram} can be thought of as a pull-back diagram of fibrations $p\times p$ and $q$ with fiber $\Omega M\times \Omega M$. 
\begin{equation*}
\begin{CD}
\Omega M\times \Omega M @= \Omega M\times \Omega M \\
@V{\iota}VV  @V{\iota'}VV \\
LM\times LM @<{g_{\text{out}}}<< \textup{Map}(\begin{tikzpicture}\draw (0,0) arc (0:180: 0.1) -- ++(- 0.2,0) arc (0:360: 0.1) -- ++( 0.2,0) arc (180:360: 0.1);\fill (- 0.2,0) circle ( 0.035); \fill (- 0.4,0) circle ( 0.035); \end{tikzpicture},M) \\
@V{p\times p}VV  @V{q}VV \\
M\times M @<{\overline{g}}<< \text{Map}(I,M)
\end{CD}
\end{equation*}
Since $H^*(LM\times LM;k)\cong H^*(M\times M;k)\otimes H^*(\Omega M\times \Omega M;k)$ by what we proved in part (1), the inclusion map of the fiber $\iota:\Omega M\times \Omega M \longrightarrow LM\times LM$ is totally nonhomologous to zero and $\iota^*$ is onto. This implies that the fiber inclusion map $\iota':\Omega M\times \Omega M \longrightarrow \textup{Map}(\begin{tikzpicture}\draw (0,0) arc (0:180: 0.1) -- ++(- 0.2,0) arc (0:360: 0.1) -- ++( 0.2,0) arc (180:360: 0.1);\fill (- 0.2,0) circle ( 0.035); \fill (- 0.4,0) circle ( 0.035); \end{tikzpicture},M)$ is such that $(\iota')^*$ is onto and $\iota'$ is totally nonhomologous to zero. Hence the spectral sequence for the fibration $q$ collapses and we get $H^*\bigl(\textup{Map}(\begin{tikzpicture}\draw (0,0) arc (0:180: 0.1) -- ++(- 0.2,0) arc (0:360: 0.1) -- ++( 0.2,0) arc (180:360: 0.1);\fill (- 0.2,0) circle ( 0.035); \fill (- 0.4,0) circle ( 0.035); \end{tikzpicture},M);k\bigr)\cong H^*(M;k)\otimes H^*(\Omega M\times \Omega M;k)$ and $g_{\text{out}}^*=\overline{g}^*\otimes 1$ is onto. By Lemma \ref{vanishing transfer}, we see that the transfer map $g_{\text{out}}^!$ vanishes. This completes the proof of Theorem \ref{cohomology of LM}.
\end{proof}

\begin{remark}\label{torsion elements and exterior algebra} If the characteristic $p$ of the field $k$ is positive, then the hypothesis \eqref{cohomology of Omega M} follows from the $p$-torsion freeness of the integral cohomology of $\Omega M$. More precisely, by Theorem 2.12 in Chapter 7 of \cite{MT}, if $H^*(\Omega M;\mathbb Z)$ is torsion free and finitely generated, then for some $\ell$, 
\begin{equation*}
H^*(\Omega M;\mathbb Z)\cong\Lambda_{\mathbb Z}(x_1,x_2,\dots,x_{\ell}), \quad |x_i|=\text{odd}. 
\end{equation*}
If $H^*(\Omega M;\mathbb Z)$ is $p$-torsion free and finitely generated, and if the characteristic of a field $k$ is $p>0$, then for some $\ell$ 
\begin{equation*}
H^*(\Omega M;k)\cong \Lambda_k(x_1,x_2,\dots,x_{\ell}), \quad |x_i|=\text{odd}.
\end{equation*}
\end{remark}

Combining the above remark and Theorem \ref{cohomology of LM}, we obtain the following corollary. 

\begin{corollary} \label{trivial coproduct} \textup{(1)} Suppose $H^*(\Omega M;\mathbb Z)$ is torsion free and finitely generated. Then the cohomology loop coproduct map in $H^*(LM;\mathbb Z)$ is trivial. 

\textup{(2)} Suppose $H^*(LM;\mathbb Z)$ has no $p$-torsion for a prime $p$. Then the cohomology loop coproduct in $H^*(LM;k)$ is trivial, where $k$ is any field of characteristic $p$. 
\end{corollary}
\begin{proof} We only have to observe that in the torsion free case (1), all the arguments in Theorem \ref{cohomology of LM} are valid with integral coefficients. 
\end{proof} 

\begin{remark}\label{differentials}
In some cases, it is not difficult to determine algebraically all the differentials in the Serre spectral sequence for the fibration $q:\text{Map}(I,M) \to M\times M$, given the cohomology of $M$. Thus the method we employ here can also be used to determine purely algebraically all the nontrivial differentials in the Serre spectral sequence for some fibrations $p:LM \to M$ without using any external geometric information. For example, all the nontrivial differentials for the fibration $p:L\mathbb CP^n \to \mathbb CP^n$ can be determined this way. In \cite{CJY},  Cohen, Jones and Yan compute the loop homology algebra structure of $H_*(LM;\mathbb Z)$ using a Serre type spectral sequence in which differentials were determined using the geometric result of Ziller \cite{Zil}, where the cohomology of $L\mathbb CP^n$ is determined using a Morse theoretic method. 
\end{remark}

\bigskip

\section{Vanishing theorem for open-closed string topology} \label{open-closed string topology}

We extend our results to open-closed string topology on an oriented closed smooth manifold $M$ of dimension $d$. First we briefly recall the general framework of open-closed string topology. An open-closed cobordism is an oriented cobordism between two compact ordered parametrized 1-dimensional manifolds, which are ordered finite union of unit intervals and unit circles. We simply refer to intervals as open strings and to circles as closed strings. Thus, an open-closed cobordism is an oriented surface $S$ whose boundary consists of three parts: (i) incoming open or closed strings $\partial_{in}S$, (ii) outgoing open or closed strings $\partial_{out}S$, and (iii) the remaining boundary part called free boundary $\partial_{free}S$. The manfold $\partial_{free}S$ is a cobordism between boundaries of $\partial_{in}S$ and $\partial_{out}S$. We assume that each connected component of an open-closed cobordism has at least one outgoing open or closed strings. This is the positive boundary condition. End points of open strings are only allowed to move along submanifolds belonging to a specified family of submanifolds of $M$ called $D$-branes. Let $\mathcal D=\{I,J,\dots\}$ be such a collection. Thus connected components of the free boundary $\partial_{free}S$ are labeled with $D$-branes. Note that a boundary component of an open-closed cobordism can have both incoming open strings and outgoing open strings, and also it can be a free boundary entirely. See \cite{Ra} and \cite{Su} for details.

The mapping class group $\Gamma(S)$ for such an open-closed cobordism $S$ is the group of isotopy classes of orientation preserving diffeomorphisms of $S$ which fix incoming and outgoing strings pointwise and which may permute completely free boundaries as long as they carry the same $D$-brane labels. Using isotopy, we may assume that such diffeomorphisms fix boundary components containing open or closed strings pointwise. 

For simplicity, suppose the set of $D$-branes consists only of $M$. For a  general case, see the remark at the end of this section. If a connected open-closed cobordism $S$ of genus $g$ has $n$ boundaries containing open or closed strings and $m$ completely free boundaries carrying the same label $M$, then the mapping class group $\Gamma(S)$ is isomorphic to $\Gamma_{g,n}^{(m)}=\pi_0(\Lambda)$, where $\Lambda$ is the topological group of orientation preserving diffeomorphisms of a genus $g$ connected oriented surface $F_{g,n}^m$ with $n$ boundaries and $m$ marked points, where diffeomorphisms fix $n$ boundaries pointwise and possibly permute $m$ marked points. Let $\Gamma_{g,n}^m$ be the mapping class group of $F_{g,n}^m$ in which diffeomorphisms fix not only $n$ boundaries pointwise but also $m$ marked points. The group $\Gamma_{g,n}^m$ is a normal subgroup of $\Gamma_{g,n}^{(m)}$ and we have the following exact sequence:
\begin{equation}\label{group extension}
1\longrightarrow \Gamma_{g,n}^m \longrightarrow \Gamma_{g,n}^{(m)} \longrightarrow \Sigma_m \longrightarrow 1.
\end{equation}
Let $\sigma_m\cong\mathbb Z$ be the sign representation of $\Sigma_m$, and we regard it as a $\Gamma_{g,n}^{(m)}$-module through the projection to $\Sigma_m$. As such $\sigma_m$ is a trivial $\Gamma_{g,n}^m$-module. Note that $\sigma_m^2$ is a trivial $\Sigma_m$-module. 

In \cite{Go2}, Godin constructs string operations in open-closed string topology in which the set of $D$-branes consists only of $M$ itself. Suppose a connected open-closed cobordism $S$ has $p$ incoming closed strings, $q$ outgoing closed strings, $r$ incoming open strings, and $s$ outgoing open strings, where we assume $q+s\ge1$, due to the positive boundary condition. Suppose the open-closed cobordism surface $S$ has genus $g$ with $n$ boundaries containing open or closed strings, and with $m$ completely free boundaries. Let $\partial_{\textup{in}}S$ be the collection of incoming open and closed strings in $S$, and let $\chi_S$ be the $\mathbb Z_2$ graded $\Gamma(S)$-module $\bigl(H_1(S,\partial_{\textup{in}}S), H_0(S,\partial_{\textup{in}}S)\bigr)$. Consider a $\Gamma(S)$ module
\begin{equation*}
\det\chi_S=\bigl(\det H_1(S,\partial_{\textup{in}}S)\bigr)\otimes \bigl(\det H_0(S,\partial_{\textup{in}}S)\bigr)^{-1}, 
\end{equation*}
where $\det$ denotes the highest exterior power. Then the associated string operations constructed in \cite{Go2} are of the following form, (modulo K\"unneth theorem):  
\begin{equation}\label{open-closed string operation}
\mu: H_*(\Gamma_{g,n}^{(m)};(\det\chi_S)^d)\otimes 
H_*(LM)^{\otimes p}\otimes H_*(M)^{\otimes r} \longrightarrow 
H_*(LM)^{\otimes q}\otimes H_*(M)^{\otimes s},
\end{equation}
where $\Gamma_{g,n}^{(m)}\cong\Gamma(S)$. 
We show that the representation $\det\chi_S$ of $\Gamma(S)$ is isomorphic to the module $\sigma_m$ described above as $\Gamma(S)$-modules. 

\begin{proposition} \label{chi_S}  Suppose and open-closed cobordism $S$ has $m$ completely free boundaries. Then as $\Gamma(S)$-module, $\det\chi_S\cong \sigma_m$. 
\end{proposition}
\begin{proof} First suppose the open-closed cobordism $S$ is connected. Then $H_0(S)\cong\mathbb Z$ is a trivial $\Gamma(S)$-module. To understand $H_1(S)$ as a $\Gamma(S)$-module, suppose $S$ has genus $g$ and has $p$ boundaries $c_1,\dots, c_p$ containing incoming open or closed strings, $m$ completely free boundaries $d_1,\dots,d_m$, and $q$ boundaries $e_1,\dots,e_q$ containing outgoing open or closed strings, where $p,m\ge0, q\ge1$ by positive boundary condition. Let $\hat{S}$ be the closed surface obtained by capping $p+m+q$ boundaries of $S$, and let $a_i,b_i$ with $1\le i\le g$ be a symplectic basis of $H_1(\hat{S})$. Then the homology basis of $H_1(S)$ consists of $a_i,b_i,[c_j],[d_k],[e_{\ell}]$, where $1\le i\le g$, $1\le j\le p$, $1\le k\le m$ and $1\le \ell\le q-1$. Note that the last homology class $[e_q]$ is a linear combination of other basis elements. Since $\Gamma(S)$ acts trivially on the basis elements $[c_j]$'s and $[e_{\ell}]$'s, and $\Gamma(S)$ permutes $[d_k]$'s, we have $\det H_1(S)\cong\bigl(\det H_1(\hat{S})\bigr)\otimes\sigma_m$ as a $\Gamma(S)$-module. Since the action of $\Gamma(S)$ on $H_1(\hat{S})$ preserves the intersection pairings of $a_i$'s and $b_j$'s, the action factors through a symplectic group, and consequently $\det H_1(\hat{S})$ is a trivial $\Gamma(S)$-module. Hence $\det H_1(S)\cong \sigma_m$ as $\Gamma(S)$-modules. 

In the general case, consider the following homology exact sequence of pairs:
\begin{equation*}
0\to H_1(\partial_{\textup{in}}S) \to H_1(S) \to H_1(S,\partial_{\textup{in}}S) 
\to H_0(\partial_{\textup{in}}S) \to H_0(S) \to  H_0(S,\partial_{\textup{in}}S) \to 0.
\end{equation*}
Since $\Gamma(S)$ acts trivially on $H_*(\partial_{\textup{in}}S)$, the above exact sequence gives 
\begin{equation*}
\det H_1(S,\partial_{\textup{in}}S)\otimes \det H_0(S)\cong \det H_1(S)\otimes \det H_0(S,\partial_{\textup{in}}S).
\end{equation*}
Let $S=\coprod_iS_i$ be the decomposition into connected components. Then $\Gamma(S)\cong\prod_i\Gamma(S_i)$ and 
\begin{equation*}
\det\chi_S\cong\det H_1(S)\otimes \bigl(\det H_0(S)\bigr)^{-1} \cong \textstyle{\bigotimes}_i\bigl(\det H_1(S_i)\otimes \det^{-1} H_0(S_i)\bigr)
\cong\textstyle{\bigotimes}_i\det H_i(S_i).
\end{equation*}
If $S_i$ has $m_i$ completely free boundaries with $m=\sum_im_i$, then by what we have proved earlier, we have $\det\chi_S\cong\bigotimes_i\sigma_{m_i}\cong\sigma_m$, as $\Gamma(S)$-modules. This completes the proof. 
\end{proof} 

If we can prove a stability property for the homology $H_*(\Gamma_{g,n}^{(m)};\sigma_m^d)$ for large genus $g$, then the same argument as in the previous section applies and we have a similar vanishing theorem for open-closed string operations. This is what we do. 

\begin{remark}\label{related stability} In \cite{CoMa} and \cite{I2}, stability properties of the homology of the mapping class group $\Gamma_{g,n}$ with twisted coefficients are studied. Here note that $\Gamma_{g,n}$ fixes boundaries of the surface $F_{g,n}$, where the above mapping class group $\Gamma_{g,n}^{(m)}$ can permute $m$ of the $m+n$ boundaries of a genus $g$ surface. So our context is somewhat different from the above papers. 
\end{remark}

A stabilizing map for an open-closed cobordism surface $S$ can be constructed in two ways. Let $T_{\text{closed}}$ be a torus with one incoming and one outgoing closed strings, and let $T_{\text{open}}$ be a torus with one boundary containing one incoming and one outgoing open strings. See Figures 2 and 3. If $S$ has an incoming or outgoing closed string, we can sew the torus $T_{\text{closed}}$ to a closed string of $S$. If $S$ has an incoming or outgoing open string, then we can sew $T_{\text{open}}$ to an open string of $S$. These two types of sewing increases the genus of $S$ by one without changing the numbers $p,q,r,s$ of incoming/outgoing open/closed strings, and the numbers $n,m$ of boundaries with or without open/closed strings. Because of the positive boundary condition $q+s\ge1$ in \eqref{open-closed string operation}, we can always apply at least one of the above two sewing procedures to increase the genus of any open-closed cobordism surface $S$.

\begin{figure}

\begin{tikzpicture}
\draw (0,0) ellipse (1.5 and 1);
\draw (0,2) node[] (1) {} arc (90:150:3  and 2 ) 
      arc (330:270:1  and 0.66 ) --++(-0.5,0) node[] (2) {};
\draw (1) arc (90:30:3  and 2 ) arc (210:270:1  and 0.66 ) 
      -- ++(0.5,0) node[] (3) {};
\draw (0,-2) node[] (4) {} arc (270:210:3  and 2 ) 
      arc (30:90:1  and 0.66 ) -- ++(-0.5,0) node[] (5) {};
\draw (4) arc (270:330:3  and 2 ) arc (150:90:1  and 0.66 ) 
      -- ++(0.5,0) node[] (6) {};
\draw[densely dashed] (0,2) arc (90:270:0.2  and 0.5 );
\draw       (0,1) arc (270:360:0.2  and 0.5 ) arc (0:90:0.2  and 0.5 );
\draw[densely dashed] (0,-1) arc (90:270:0.2  and 0.5 );
\draw    (0,-2) arc (270:360:0.2  and 0.5 ) arc (0:90:0.2  and 0.5 );
\draw[ultra thick] (-3.95,0) ellipse (0.25  and 0.65 );
\draw[ultra thick] (3.95,0) ellipse (0.25  and 0.65 );
\path (0,-3) node[text width=10cm] 
{\textsc{Figure 2.} The torus $T_{\textup{closed}}$ has one incoming boundary and one outgoing boundary.};

\end{tikzpicture}

\qquad
\begin{tikzpicture}[>=stealth]
\draw (2.5,0) arc (0:60: 2.5  and 1 );
\draw (-2.5,0) arc (180:120:2.5  and 1 );
\draw (-2.5,0) arc (180:360:2.5  and 1 );
\draw[ultra thick] (2.5,0) arc (0:30:2.5  and 1 );
\draw[ultra thick] (-2.5,0) arc (180:210: 2.5  and 1 );
\draw[ultra thick] (-2.5,0) arc (180:150: 2.5  and 1 );
\draw[ultra thick] (2.5,0) arc (360:330:2.5  and 1 );
\draw (-1.5,0) arc (180:360:0.5  and 0.2 ) 
arc (180:0:0.5 ) arc (180:360:0.5  and 0.2 ) 
arc (0:180:1.5 );
\draw[->, ultra thick] (2.5,0) -- ++(0,0.1);
\draw[->, ultra thick] (-2.5,0) -- ++(0,0.1);
\path (1.8,-0.7) node[below] {$I$};
\path (1.8,0.7) node[above] {$J$};
\draw[dashed] (1.5,0) arc(0:180:0.5  and 0.2 );
\draw[dashed] (-0.5,0) arc (0:180:0.5  and 0.2 );
\node[text width=10cm] at (0,-1.7) 
{\textsc{Figure 3}. The torus $T_{\textup{open}}$ has one boundary with one incoming open string and one outgoing open string};
\end{tikzpicture}

\end{figure} 

By extending diffeomorphisms on $S$ by identity on $T_{\text{closed}}$ or on $T_{\text{open}}$, we obtain a homomorphism $\varphi: \Gamma(S)\rightarrow \Gamma(S\#T)$, that is, a homomorphism $\varphi:\Gamma_{g,n}^{(m)} \to \Gamma_{g+1,n}^{(m)}$. By restriction to those diffeomorphisms of $S$ preserving completely free boundaries component wise, we obtain a homomorphism $\varphi: \Gamma_{g,n}^m \to \Gamma_{g+1,n}^m$. 
We show that Harer-Ivanov's stability result in the introduction can be extended to the above mapping class groups with the same stability range and with the coefficient in the module $(\sigma_m)^{\otimes r}$. 

For this, we first recall Ivanov's formulation of stability theorem. Let $X,Y$ be compact connected orientable surfaces with nonempty boundary such that $X\subset Y$. Let $\mathcal M_X$ and $\mathcal M_Y$ be their mapping class group consisting of isotopy classes of orientation preserving diffeomorphisms fixing the boundaries pointwise. Thus, if $X$ is a surface of genus $g$ with $n$ boundaries, then $\mathcal M_X\cong\Gamma_{g,n}$. Let $g(X)$ be the genus of the surface $X$. 

\begin{theorem}[Ivanov \cite{I2}]\label{Ivanov stability} With the above notation, the homomorphism of mapping class groups induced by the inclusion $X\longrightarrow Y$
\begin{equation*}
H_k\bigl(\mathcal M_X\bigr)  \longrightarrow H_k\bigl(\mathcal M_Y\bigr)
\end{equation*}
is an isomorphism when $g(X)\ge 2k+1$, and onto when $g(X)\ge 2k$. 
\end{theorem}

In the original formulation of stability theorem by Harer, it uses sewing of the torus $T_{\text{closed}}$ along a boundary of a surface, and it does not cover the case of sewing $T_{\text{open}}$ to a surface. This latter case is covered by Ivanov's formulation of the stability theorem. When the surface $Y$ is obtained by sewing $T_{\text{closed}}$ or $T_{\text{open}}$ to $X$, the actual homomorphism may depend on the choice of (part of) the boundary of $X$ used for sewing, as we saw in the previous section.

\begin{theorem} \label{homology stability} Let $n\ge1$ and $r\ge0$. Let $\varphi$ be a stabilizing map obtained by sewing $T_{\textup{closed}}$ or $T_{\textup{open}}$ to the open-closed cobordism $S$. Then, both of the following homology stabilizing maps are isomorphisms for $g\ge 2k+1$ and are onto for $g\ge2k$\textup{:}
\begin{align}
\varphi_*&: H_k(\Gamma_{g,n}^m) \longrightarrow H_k(\Gamma_{g+1,n}^m), 
\label{stability with marked points}\\
\varphi_*&: H_k(\Gamma_{g,n}^{(m)};\sigma_m^r) \longrightarrow H_k(\Gamma_{g+1,n}^{(m)};\sigma_m^r).\label{stability with twisted module}
\end{align}
Here $\sigma_m^r=(\sigma_m)^{\otimes r}$ is a trivial $\Gamma_{g,n}^{(m)}$-module for even $r$. 

When $g\ge 2k+1$, the homology groups $H_k(\Gamma_{g,n}^m)$ and $H_k(\Gamma_{g,n}^{(m)})$ are independent of $n\ge1$. Furthermore, when $g\ge2k$, the action of $\Sigma_n$ on both homology groups $H_k(\Gamma_{g,n}^m)$ and $H_k(\Gamma_{g,n}^{(m)})$ is trivial and consequently, the stabilizing map $\varphi_*$ in \eqref{stability with marked points} and \eqref{stability with twisted module} is independent of the choice of boundaries of $S$ used for sewing with $T_{\textup{closed}}$ or $T_{\textup{open}}$. 
\end{theorem}
\begin{proof} For the first case, we have the following exact sequence (see for example \cite{ABE}):
\begin{equation}\label{group extension: capping by discs}
1 \longrightarrow \mathbb Z^m \longrightarrow \Gamma_{g,n+m} \longrightarrow \Gamma_{g,n}^m  \longrightarrow 1, 
\end{equation}
where $\Gamma_{g,n+m}$ is the mapping class group of $S$ consisting of isotopy classes of  orientation preserving diffeomorphisms of $S$ fixing all $n+m$ boundaries pointwise. This is a central extension and the kernel $\mathbb Z^m$ is generated by Dehn twists along simple closed curves parallel to $m$ completely free boundaries. We consider the associated Hochschild-Serre spectral sequence
\begin{equation*}
E^2_{p,q}=H_p\bigl(\Gamma_{g,n}^m;H_q(\mathbb Z^m)\bigr) \Longrightarrow 
H_{p+q}(\Gamma_{g,n+m}).
\end{equation*}
Since the above extension is a central extension, the action of $\Gamma_{g,n}^m$ on $\mathbb Z^m$ is trivial, and thus we have a trivial local system in the above $E^2$-terms. Since the homology $H_*(\mathbb Z^m)$ is torsion free, the above $E^2$-term can be written as $E^2_{p,q}=H_p(\Gamma_{g,n}^m)\otimes H_q(\mathbb Z^m)$. Now sewing $T_{\text{closed}}$ or $T_{\text{open}}$ to $S$ induces the following homomorphisms between group extensions.
\begin{equation*}
\begin{CD}
1 @>>> \mathbb Z^m @>>> \Gamma_{g,n+m} @>>> \Gamma_{g,n}^m @>>> 1 \\
 @. @| @V{\varphi}VV @V{\varphi}VV @. \\
1 @>>> \mathbb Z^m @>>> \Gamma_{g+1,n+m} @>>> \Gamma_{g+1,n}^m @>>> 1
\end{CD}
\end{equation*}
This diagram induces a homomorphism of spectral sequences 
\begin{equation*}
E^2_{p,q}=H_p(\Gamma_{g,n}^m)\otimes H_q(\mathbb Z^m) 
\xrightarrow{\varphi_*\otimes 1} {}'E^2_{p,q}=H_p(\Gamma_{g+1,n}^m)\otimes 
H_q(\mathbb Z^m)
\end{equation*}
converging to the homomorphism $\varphi_*: H_{p+q}(\Gamma_{g,n+m}) \longrightarrow H_{p+q}(\Gamma_{g+1,n+m})$ which we know to be an isomorphism for $g\ge 2(p+q)+1$ and onto for $g\ge2(p+q)$ by Harer-Ivanov stability theorem. Using one version of Zeeman's comparison theorem of spectral sequences (see Theorem 1.3 in \cite{I2}), stability property of the group $\Gamma_{g,n+m}$ (via sewing of $T_{\text{closed}}$ or $T_{\text{open}}$ to open-closed cobordisms $S$) implies the stability property of $\Gamma_{g,n}^m$ in the same range. Namely, $\varphi_*:H_k(\Gamma_{g,n}^m) \rightarrow  H_k(\Gamma_{g+1,n}^m)$ is an isomorphism for $g\ge 2k+1$ and onto for $g\ge 2k$. This proves the first part. 

For the second stabilizing map, we again consider Hochschild-Serre spectral sequence associated to the group extension \eqref{group extension} and the $\Gamma_{g,n}^{(m)}$-module $\sigma_m^r$:  
\begin{equation*}
E^2_{p,q}=H_p\bigl(\Sigma_m; H_q(\Gamma_{g,n}^m;\sigma_m^r)\bigr) \Longrightarrow H_{p+q}(\Gamma_{g,n}^{(m)};\sigma_m^r),  
\end{equation*}
where $H_*(\Gamma_{g,n}^m;\sigma_m^r)=H_*(\Gamma_{g,n}^m;\mathbb Z)\otimes\sigma_m^r$ since $\Gamma_{g,n}^m$ acts trivially on the module $\sigma_m^r$. Now, sewing $T_{\text{closed}}$ or $T_{\text{open}}$ to the open-closed cobordism $S$ induces the following homomorphism between group extensions: 
\begin{equation}\label{extension diagram}
\begin{CD}
1 @>>> \Gamma_{g,n}^m @>>> \Gamma_{g,n}^{(m)} @>>> \Sigma_m @>>> 1  \\
@. @V{\varphi}VV @V{\varphi}VV @|  @.\\
1 @>>> \Gamma_{g+1,n}^m @>>> \Gamma_{g+1,n}^{(m)} @>>> \Sigma_m @>>> 1, 
\end{CD}
\end{equation}
which induces a homomorphism of spectral sequences 
\begin{equation*}
E^2_{p,q}=H_p\bigl(\Sigma_m; H_q(\Gamma_{g,n}^m)\otimes\sigma_m^r\bigr) \longrightarrow 
{}'E^2_{p,q}=H_p\bigl(\Sigma_m; H_q(\Gamma_{g+1,n}^m)\otimes\sigma_m^r\bigr)
\end{equation*}
converging to the homomorphism $\varphi_*: 
H _{*}(\Gamma_{g,n}^{(m)};\sigma_m^r) \longrightarrow 
H_{*}(\Gamma_{g+1,n}^{(m)};\sigma_m^r)$. A standard Zeeman's comparison theorem of spectral sequence (see Theorem 1.2 in \cite{I2}) together with the stability property for $\Gamma_{g,n}^m$ we have just proved, we conclude that the group $\Gamma_{g,n}^{(m)}$ also enjoys a stability property in the stated range. 

To see that the homology group $H_k(\Gamma_{g,n}^m)$ is independent of $n\ge1$ when $g\ge 2k+1$, we consider the following diagram:
\begin{equation*}
\begin{CD}
1 @>>> \mathbb Z^m @>>> \Gamma_{g,n+m} @>>> \Gamma_{g,n}^m @>>> 1 \\
@. @| @VVV @VVV @. \\
1 @>>> \mathbb Z^m @>>> \Gamma_{g,1+m} @>>> \Gamma_{g,1}^m @>>> 1, 
\end{CD}
\end{equation*}
where the vertical maps are induced by capping $n-1$ incoming boundaries with discs. By Theorem \ref{Ivanov stability}, the induced middle vertical homomorphisms in homology is onto when $g\ge 2k$ and isomorphism when $g\ge 2k+1$. Thus, using Zeeman's spectral sequence comparison theorem as above, we see that the same is true for induced homology map on the right. Thus, for $g\ge 2k+1$, we have $H_k(\Gamma_{g,1}^m)\cong H_k(\Gamma_{g,n}^m)$ for any $n\ge1$. 

In the above context, if we sew a surface $F_{0,1+n}$ to $F_{g,1+m}$ or to $F_{g,1}^m$, then we have homomorphisms going from the bottom row to the top row. We can deduce the same conclusion arguing as before using this new diagram with reversed vertical arrows. 

The proof that the homology group $H_k(\Gamma_{g,n}^{(m)})$ is independent of $n\ge1$ when $g\ge2k+1$ is the same as above using analogous diagrams induced by either capping $n-1$ incoming boundaries or sewing $F_{0,1+n}$ to $F_{g,1+m}$ and to $F_{g,1}^{(m)}$. 

Finally, triviality of the action of $\Sigma_n$ on homology groups $H_k(\Gamma_{g,n}^m)$ and $H_k(\Gamma_{g,n}^{(m)})$ can be shown in exactly the same way as in Remark \ref{trivial S_n action} in the stable range $g\ge2k$. Then arguing as in Proposition \ref{transposition}, we can see that the stabilizing map $\varphi_*$ in \eqref{stability with marked points} and \eqref{stability with twisted module} is independent of choices involved in sewing the open-closed cobordism $S$ and $T_{\text{closed}}$ or $T_{\text{open}}$ in the same stable range $g\ge2k$. This completes the proof. 
\end{proof}

In the Harer's original paper \cite{H}, stability property of the group $\Gamma_{g,n}^m$ is proved, but the slope of the stability range is $3$, instead of $2$ as in Theorem \ref{homology stability}. 

The proof of the stability property of the stabilizing map \eqref{stability with twisted module} was done in two steps, using a spectral sequence comparison theorem each time. We could complete the proof of the stability of \eqref{stability with twisted module} using a single spectral sequence. Let $F_{g,n}^{(m)}$ be a smooth oriented surface of genus $g$ with $n$ boundaries, $m$ marked points $\{x_1,x_2,\dots,x_m\}$, and a choice of an oriented frame $(u_i,v_i)$ at each marked point $x_i$. Let $\text{Diff}^+\bigl(F_{g,n}^{(m)}\bigr)$ be the topological group of orientation preserving diffeomorphisms which fix $n$ boundaries pointwise and which possibly permute $m$ marked points. For each diffeomorphism $f\in \text{Diff}^+\bigl(F_{g,n}^{(m)}\bigr)$, let $f(x_i)=x_{\tau(i)}$, $1\le i\le m$, for some permutation $\tau\in\Sigma_m$, and let $\bigl(f_*(u_i),f_*(v_i)\bigr)=(u_{\tau(i)},v_{\tau(i)})A_i$ for some $A_i\in\text{GL}^+_2(\mathbb R)$ for $1\le i\le m$. This correspondence from $f$ to $(\tau;A_1,A_2,\dots,A_m)$ defines an onto homomorphism from the diffeomorphism group to a wreath product $\text{Diff}^+\bigl(F_{g,n}^{(m)}\bigr) \longrightarrow \Sigma_m\wr \text{GL}^+_2(\mathbb R)$ whose kernel consists of those diffeomorphisms which fix $m$ marked points and whose induced differentials are identity at $m$ marked points. Thus the kernel is homotopy equivalent to $\text{Diff}^+(F_{g,n+m})$. Since $n\ge1$, connected components of these diffeomorphism groups are contractible \cite{ES}. Thus, passing to classifying spaces and then replacing diffeomorphism groups by mapping class groups, we have a homotopy fibration given in the top row of the next diagram, and a map of fibrations induced by a stabilizing map as in the diagram below:
\begin{equation*}
\begin{CD}
B\Gamma_{g,n+m} @>>>  B\Gamma_{g,n}^{(m)} @>>> B\bigl(\Sigma_m\wr\text{GL}^+_2(\mathbb R)\bigr) \\
@V{B\varphi}VV  @V{B\varphi}VV @| \\
B\Gamma_{g+1,n+m} @>>>  B\Gamma_{g+1,n}^{(m)} @>>> B\bigl(\Sigma_m\wr\text{GL}^+_2(\mathbb R)\bigr).
\end{CD}
\end{equation*}
Here, since $\text{GL}^+_2(\mathbb R)$ is homotopy equivalent to the circle $S^1$, the classifying space $B\bigl(\Sigma_m\wr\text{GL}^+_2(\mathbb R)\bigr)$ is homotopy equivalent to $E\Sigma_m\times_{\Sigma_m}(\mathbb CP^{\infty})^m$. Now we consider maps between two Serre spectral sequences associated to the above two fibrations in the top and the bottom rows, and by applying Zeeman's comparison theorem, we get the same result as in Theorem \ref{homology stability}. See \cite{BT} for related results on stable mapping class groups. 

As before, we say that the group $H_k(\Gamma_{g,n}^{(m)};\sigma_m^d)$ is in stable range if the stabilizing map $\varphi_*$ mapping to this group is onto and all the subsequent stabilizing maps are isomorphisms: 
\begin{equation*}
H_k(\Gamma_{g-1,n}^{(m)};\sigma_m^d) \xrightarrow[\text{onto}]{\varphi_*} 
H_k(\Gamma_{g,n}^{(m)};\sigma_m^d)
\xrightarrow[\cong]{\varphi_*} 
H_k(\Gamma_{g+1,n}^{(m)};\sigma_m^d)
\xrightarrow[\cong]{\varphi_*} \cdots 
\end{equation*}

\begin{Vanishing Theorem}[\textbf{Open-Closed String Topology Case}] Let $M$ be a smooth oriented closed manifold of dimension $d$. Consider open-closed string topology on $M$ in which the set of $D$-brane submanifolds consists only of $M$. Let 
\begin{equation*}
\varphi_*: 
H_k(\Gamma_{g,n}^{(m)};\sigma_m^d) \longrightarrow 
H_k(\Gamma_{g+1,n}^{(m)};\sigma_m^d)
\end{equation*}
be a stabilizing map of mapping class groups obtained by sewing either $T_{\textup{closed}}$ or $T_{\textup{open}}$ to open-closed cobordisms. 

\textup{(i)} Open-closed string operations \eqref{open-closed string operation} associated to elements in the image $\textup{Im}\,\varphi_*$ of any stabilizing map $\varphi_*$ are trivial. 

\textup{(ii)} Open-closed string operations \eqref{open-closed string operation} associated to any elements in the homology group  $H_k(\Gamma_{g,n}^{(m)};\sigma_m^d)$ in stable range are trivial. 
\end{Vanishing Theorem}
\begin{proof} As before, using the gluing property of the homological conformal field theory, we only have to observe that the topological quantum field theory operations associated to surfaces $T_{\textup{closed}}$ and $T_{\textup{open}}$ are trivial.  We have already seen that the TQFT operation associated to $T_{\textup{closed}}$ is trivial. The TQFT operation associated to $T_{\textup{open}}$ describs a process in which a closed string splits off from an open string, and later they join together to form an open string. The general form of such string operation in which open strings carry submanifold labels $I,J$ at their end points is given by
\begin{equation*}
\mu_{T_{\textup{open}}}: H_*(P_{IJ}) \longrightarrow H_*(P_{IJ})\otimes H_*(LM)  \longrightarrow H_*(P_{IJ}),
\end{equation*}
where $P_{IJ}$ is the space of open strings $\gamma:[0,1] \rightarrow M$ such that $\gamma(0)\in I$ and $\gamma(1)\in J$. In Proposition 3.4 in \cite{T5} we called such an operation a handle attaching operation, and we showed that handle attaching operations are always trivial for any labels $I,J$. In our present case, we have $I=J=M$. By a similar argument as in the closed string topology case, the open-closed string operation associated to any element in the image of any stabilizing map $\varphi_*$ can be factored into a composition of operations, and one of the factors is the TQFT operation associated to $T_{\textup{closed}}$ or $T_{\textup{open}}$, which is trivial. This proves part (i). Part (ii) is a consequence of part (i) and definition of stable range. 
\end{proof}  

As a final remark in this section, we consider a general case in which the set $\mathcal D$ of $D$-brane submanifolds has more than one element. Let $\mathcal D=\{K_1,K_2,\dots, K_h,\dots\}$ be a set of $D$-brane labels. Let $S_{g,n}^m$ be an oriented open-closed cobordism of genus $g$ with $n$ boundaries containing open or closed strings, and with $m$ completely free boundaries. In this case, we consider orientation preserving diffeomorphisms of $S_{g,n}^m$ which fix $n$ boundaries containing open or closed strings pointwise, and which may permute completely free boundaries provided they carry the same label. Suppose there are $m_i$ completely free boundaries carrying the same label $K_i$ for $1\le i\le h$ so that $\sum_im_i=m$. Let $H(\vec{m})=\Sigma_{m_1}\times \Sigma_{m_2}\times\cdots\times \Sigma_{m_h}\subset\Sigma_m$ be the subgroup of the symmetric group corresponding to $\vec{m}=(m_1,m_2,\dots,m_h)$. The corresponding mapping class group is denoted by $\Gamma_{g,n}^{H(\vec{m})}$. Then as before, there exists a homotopy fibration 
\begin{equation*}
B\Gamma_{g,n+m} \longrightarrow B\Gamma_{g,n}^{H(\vec{m})} \longrightarrow B\bigl(H(\vec{m})\wr\text{GL}^+_2(\mathbb R)\bigr).
\end{equation*}
Then arguing as before, we can show that the homology group $H_k\bigl(\Gamma_{g,n}^{H(\vec{m})};\sigma_m^r\bigr)$ for $r\ge0$ has a stability property with respect to the genus $g$ as in Theorem \ref{homology stability} with the same range of stability. This stability property would be relevant to a  vanishing property of open-closed string operations with a general set of $D$-branes.

\bigskip

\section{Some computations of unstable closed string operations}

We have shown that stable string operations vanish. We can ask: what about unstable string operations? Part (i) of vanishing theorems applies to stable as well as unstable string operations, and it shows that most of the unstable string operations vanish. Those unstable operations not covered by part (i) are those operations associated to homology elements of mapping class groups not in the image of any stabilizing maps. Can they be nontrivial? We can try to compute some unstable string operations. Unfortunately not many homology groups of mapping class groups in unstable range have been calculated (see for example \cite{ABE}), although the situation for stable mapping class groups is much better (\cite{Ga} and \cite{MW}). 

In this final section, we compute some genus one unstable string operations in closed string topology for finite dimensional $M$. First we examine string operations associated to degree $1$ homology group $H_1(\Gamma_{1,p+q})$ of genus $1$ mapping class groups with $p+q\ge1$. This homology group is well known to be $H_1(B\Gamma_{1,p+q})\cong\mathbb Z^{p+q}$, generated by a Dehn twist along a nonseparating simple closed curve on the surface $F_{1,p+q}$ (we can use the meridian of the torus), and Dehn twists along simple closed curves parallel to $p+q-1$ boundary circles (for an explanation, see for example \cite{K}, Theorem 5.1). Note that it is well known that Dehn twists along nonseparating simple closed curves on any connected surface are conjugate to each other in its mapping class group, and hence they represent the same first homology classes. The string operation corresponding to the Dehn twist on a cylinder generating $H_1(\Gamma_{0,1+1})\cong\mathbb Z$ is the BV operator $\Delta:H_*(LM) \to H_{*+1}(LM)$ coming from the homological circle action on the free loop space $LM$. By inserting BV operator appropriately on a decomposition of the genus 1 surface $F_{1,p+q}$, and using the gluing property of HCFT, we see that all string operations associated to $H_1(\Gamma_{1,p+q})$ vanish. For example, consider the Dehn twist along the meridian of $F_{1,p+q}$. Let $\Psi$ and $\mu$ be the loop coproduct and the loop product maps both of degree $-d$:
\begin{align*}
\Psi&: H_*(LM) \longrightarrow H_*(LM)\otimes H_*(LM), \\
\mu &: H_*(LM)\otimes H_*(LM) \longrightarrow H_*(LM).
\end{align*}
In \cite{T3} Theorem 2.5, we showed that the coproduct is nontrivial only on $H_d(LM)$ and its image under $\Psi$ are integral multiples of the generator $[c_0]\otimes [c_0]\in H_0(LM)\otimes H_0(LM)$ by degree reason, where $c_0$ is a constant loop. The string operation associated to the Dehn twist along the meridian is given by 
\begin{equation*}
\mu\circ(\Delta\otimes 1)\circ\Psi:H_*(LM) \longrightarrow H_{*-d+1}(LM), 
\end{equation*}
and since $\Delta([c_0])=0$, this string operation is trivial. For Dehn twists along a curve parallel to boundaries, the situation is even more straightforward and they are given by post or pre-composition of the genus one operator $\mu\circ\Psi$ with a BV operator $\Delta$, and hence they are trivial, too. 

Next we compute a string operation associated to a degree $2$ homology class. Since the torus $T_{\text{closed}}$ is used in the stabilization map in Harer's theorem, we examine string operations associated to the homology group $H_*(\Gamma_{1,1+1})$. This homology group is computed in \cite{Go1}, and is given by $H_1\cong\mathbb Z\oplus\mathbb Z$, $H_2\cong\mathbb Z\oplus\mathbb Z_2$, $H_3\cong\mathbb Z_2$, and $H_k=0$ for $k\ge4$. As before the string operations associated to $H_1$ vanish. We examine the string operation associated to a generator the infinite cyclic group in $H_2$. To understand this generator, we compute the homology of $\Gamma_{1,2}$ in a different way. We have the following group extension
\begin{equation}\label{group extension: one cap}
1\longrightarrow \mathbb Z \longrightarrow \Gamma_{1,2} \longrightarrow \Gamma_{1,1}^1 \longrightarrow 1,
\end{equation}
where the homomorphism $\Gamma_{1,2} \rightarrow \Gamma_{1,1}^1$ is induced by capping one boundary of $T$ with a disc keeping the center point of the disc. The kernel $\mathbb Z$ is generated by the Dehn twist along a simple closed curve parallel to the capped boundary, and it is in the center of $\Gamma_{1,2}$. This group extension is a special case of \eqref{group extension: capping by discs} Since the action of $\Gamma_{1,1}^1$ on $\mathbb Z$ is trivial, the local system in the following Hochschild-Serre spectral sequence is trivial: 
\begin{equation*}
E^2_{p,q}=H_p\bigl(\Gamma_{1,1}^1; H_q(\mathbb Z)\bigr) \Longrightarrow H_{p+q}(\Gamma_{1,2}).
\end{equation*}
Homology of $\Gamma_{1,1}^1$ is listed in the table in \cite{ABE} as follows: $H_1(\Gamma_{1,1}^1)\cong\mathbb Z$ generated by a Dehn twist along any nonseprating simple closed curve (we can take this to be a meridian) on the genus $1$ surface with one boundary and one puncture, and $H_2(\Gamma_{1,1}^1)\cong\mathbb Z_2$. By glancing at the $E^2$-terms of the above spectral sequence, we see that there cannot be any nontrivial differentials by degree reason, and the spectral sequence must collapse. For the extension problem, the group $H_2(\Gamma_{1,2})$ fits into the following exact sequence:  
\begin{equation*}
0\longrightarrow E^{\infty}_{1,1} \longrightarrow H_2(\Gamma_{1,2}) \longrightarrow E^{\infty}_{2,0} \longrightarrow 0, 
\end{equation*}
where $E^{\infty}_{1,1}\cong\mathbb Z$ and $E^{\infty}_{2,0}\cong\mathbb Z_2$. From this exact sequence, we see that a generator of the infinite cyclic group in $H_2(\Gamma_{1,2})\cong\mathbb Z\oplus\mathbb Z_2$ comes from a generator of $E^2_{1,1}=H_1(\Gamma_{1,1}^1)\otimes H_1(\mathbb Z)\cong\mathbb Z$. Thus in the fibration $S^1\to B\Gamma_{1,2} \to B\Gamma_{1,1}^1$ associated to \eqref{group extension: one cap}, an infinite cyclic generator of $H_2(B\Gamma_{1,2})$ is given by a cycle $S^1\times S^1 \to B\Gamma_{1,2}$ where the fist $S^1$ corresponds to the Dehn twist along the meridian of the torus $T_{\text{closed}}$ corresponding to a generator of $H_1(B\Gamma_{1,1}^1)$, and the second $S^1$ corresponds to the Dehn twist along a simple closed curve parallel to one of the boundaries of $T_{\text{closed}}$, corresponding to a generator of $H_1(B\mathbb Z)=H_1(S^1)$. Thus the corresponding string operation is given by  
\begin{equation*}
\Delta\circ\mu\circ(\Delta\otimes 1)\circ\Psi :H_*(LM) \longrightarrow H_{*-d+2}(LM). 
\end{equation*}
Since as before the coproduct $\Psi$ in $H_*(LM)$ followed by a BV operator $\Delta\otimes 1$ vanishes, the string operation associated to the generator of $\mathbb Z\subset H_2(\Gamma_{1,2})\cong\mathbb Z\oplus\mathbb Z_2$ is trivial. 

For the other genus $1$ surface $T_{\text{open}}=F_{1,1}$ we used for stabilizing maps, the homology of the corresponding mapping class group $\Gamma_{1,1}$ is given by $H_1(\Gamma_{1,1})\cong\mathbb Z$ and $H_k(\Gamma_{1,1})=0$ for $k\ge2$ \cite{ABE}. Thus, there are no interesting homology classes in this case. 

For closed string topology for simply connected infinite dimensional manifold $M$ with finite dimensional $H^*(\Omega M;k)$, the discussion goes through in parallel, and we obtain the same result. 

We record our results of the above discussion on unstable genus one closed string operations in the next proposition. Of course similar results can be obtained in the context of unstable genus one open-closed string operations using similar decompositions of open-closed cobordisms. 

\begin{proposition} Let $\Gamma_{1,r}$ be the mapping class group of genus $1$ surface with $r$ boundaries. Then in closed string topology for finite or infinite dimensional $M$, the followings hold\textup{:}

\textup{(1)} Closed string operation associated to an arbitrary element in the first homology group $H_1(B\Gamma_{1,r})\cong\mathbb Z^r$ for $r\ge1$ vanishes. 

\textup{(2)} Closed string operation associated to a generator of the free summand of the second homology group $H_2(B\Gamma_{1,2})\cong\mathbb Z\oplus\mathbb Z_2$ vanishes. 
\end{proposition}

\bibliography{bibliography}
\bibliographystyle{plain}

\end{document}